\documentclass[]{siam/siamonline250211}

\usepackage{setspace}
\usepackage{bmk_math}
\usepackage{bmk_plots}
\usepackage[caption=false]{subfig}

\usepackage{amsfonts}
\usepackage{graphicx}
\usepackage{epstopdf}
\usepackage{algorithmic}
\ifpdf
  \DeclareGraphicsExtensions{.eps,.pdf,.png,.jpg}
\else
  \DeclareGraphicsExtensions{.eps}
\fi

\usepackage{enumitem}
\setlist[enumerate]{leftmargin=.5in}
\setlist[itemize]{leftmargin=.5in}

\newsiamremark{remark}{Remark}
\newsiamremark{hypothesis}{Hypothesis}
\crefname{hypothesis}{Hypothesis}{Hypotheses}
\newsiamthm{claim}{Claim}
\newsiamremark{fact}{Fact}
\crefname{fact}{Fact}{Facts}

\headers{Learning calibrated probabilistic forecasts using KMD}{B. M. Kent and R. E. Turner}

\title{Learning calibrated probabilistic forecasts using Koopman mode decomposition and the continuous ranked probability score \thanks{
\funding{This work was funded by \bk{EPSRC Probabilistic AI Hub (EP/Y028783/1).}}}}

\author{Benjamin M. Kent\thanks{University of Cambridge 
  (\email{bmk28@cam.ac.uk}, \url{benmkent.github.io}).}
\and Richard E. Turner\thanks{University of Cambridge
  (\email{ret26@cam.ac.uk}).}}

\usepackage{amsopn}

\date{}
\usepackage{hyperref}
\usepackage{cleveref}
\AddToHook{env/remark/begin}{\crefalias{theorem}{remark}}
\usepackage{placeins}
\usepackage{booktabs}
\usepackage{tabularx}
\usepackage{pgfplotstable}
\usepackage{multirow}
\usepackage{rotating}
\usepackage{standalone}
\usepackage{microtype}

\newcommand{\bk}[1]{{#1}}
\newcommand{\bmk}[1]{{#1}}

\begin{document}
\maketitle
\begin{abstract}
\bk{Well-calibrated forecasting for non-linear dynamical systems is essential in science and engineering applications such as fluid dynamics and weather prediction.
We present a forecasting method combining probabilistic machine learning with the structured linear dynamics of Koopman mode decomposition.}
\bk{Within this framework, the} observation, linear dynamics, and state reconstruction operators are constructed using standard machine learning \bk{architectures,} producing predictive distributions via the pushforward of \bk{injected} noise samples.
The model is trained using a proper scoring rule, the continuous ranked probability score, to directly target calibrated, accurate forecasts.
The approach is validated on two non-linear {ordinary differential equation} examples and two {partial differential equation} examples, vortex shedding and cavity flow, \bk{using both large and small, noisy datasets.}
\bk{Results} demonstrate that known Koopman eigenfunctions and \bk{spectra} can be identified and used for accurate forecasting.
The model can also produce multi-modal predictive distributions, demonstrated for \bk{a trajectory initialised on the boundary between two} basins of attraction.
\bk{Uncertainty} is captured by the sample ensemble spread, with spread--skill ratios close to one over the training \bk{horizon in several experiments,} although for extended \bk{rollouts and noisy datasets,} the ensemble may become under- or \bk{over-dispersed.}
\bk{In low-noise experiments, temporal rollout is generally accurate for at least twice the training horizon, and up to twenty times in some cases.}
\bk{Approximation quality is comparable to a deep ensemble approximation over the training horizon, but at a reduced computational cost and with better robustness in identifying Koopman eigenvalues.}
\end{abstract}



\section{Introduction}
\label{sec:introduction}

Prediction and forecasting for dynamical systems are essential tools in many modern scientific applications, with a growing emphasis on data-driven approximation methods~\cite{Mishra2024,Bramburger2024,BruntonKutz2022}.
For example, one may wish to understand the dynamics of vortex induced vibration when designing a structure~\cite{Raissi2018}, analyse cavity flow in industrial and scientific contexts~\cite{Shankar2000}, or forecast the weather~\cite{Alet2025,Allen2025}. 
To make well-informed decisions, it is essential to not only estimate the future state, but \bk{also} to provide a well-calibrated forecast quantifying \bk{uncertainty, and an interpretable model to aid system identification.} 
\bmk{Modelling uncertainty} is particularly important \bmk{in} time-dependent problems, where uncertainty and approximation error can accumulate rapidly over temporal rollouts.
In the field of machine learning, probabilistic modelling techniques are commonly used for prediction and forecasting under uncertainty~\cite{Price2025,Allen2025,Lang2026,diaconu2026probabilisticretrofittinglearnedsimulators,Gal2016,Lakshminarayanan2017}.
\bmk{We desire an interpretable, stable and uncertainty-aware approximation, so we present a data-driven, probabilistic approach that builds upon scoring rule based optimisation~\cite{Gneiting2007} to target calibration, and exploits the mathematical structure of Koopman-related methods~\cite{Brunton2022,Colbrook2024}.} %

The \bk{Koopman operator is a linear operator describing the evolution of observables of a dynamical system~\cite{Koopman1931,Koopman1932,Mezic2005}.
Its numerical approximation is related to Dynamic Mode Decomposition (DMD)~\cite{Rowley2009,Schmid2010,Tu2014}.}
\bmk{DMD} defines a data-based numerical method that identifies coherent state modes and corresponding linear time evolution via a least-squares error fit on trajectory data.
\bmk{DMD pairs linear dynamics} with linear projection into a \bk{truncated} state basis to lower computational cost and is related to reduced order modelling techniques such as proper orthogonal decomposition~\cite{Berkooz1993}.
Variations on DMD offer verification~\cite{Colbrook2023res}, uncertainty quantification via an ensemble of models~\cite{Sashidhar2022} and probabilistic formulations~\cite{Takeishi2017,Morton2019,Pan2020}, plus variants for control~\cite{Proctor2016,Korda2018}, parametric problems~\cite{Andreuzzi2023,Huhn2023}, and \bk{more~\cite{Askham2018,Baddoo2023,Brunton2015,Brunton2022,Colbrook2024,Demo2018}.}
\emph{Extended dynamic mode decomposition (EDMD)}~\cite{Williams2015} instead identifies linear dynamics via a least-squares error fit for \bk{a lifting of the state} into an observation space defined via a dictionary of non-linear observable functions mapping the dynamical system state to complex values in $\C$.
Under certain conditions, the EDMD approximation converges to the Koopman operator when taking two successive limits: infinite state \bk{trajectory} data followed by infinite \bk{dictionary} dimension~\cite{Korda2017}.
For dynamical systems with high-dimensional state spaces, choosing an appropriate dictionary is tricky, \bk{which can be partially addressed via kernel-based approaches~\cite{Williams2015b}}.
\bk{Alternatively,} \bk{machine} learning approaches identify a suitable observation space and corresponding linear dynamics simultaneously, generally using a least-squares error loss~\cite{Li2017,Lusch2018,Yeung2019,Takeishi2017b,Erichson2019,Baddoo2022,Otto2019}.
A Bayesian DMD~\cite{Takeishi2017} formulation builds upon probabilistic \bk{PCA~\cite{Tipping1999} and probabilistic} methods for non-linear dictionaries are investigated \bk{in~\cite{Morton2019,Pan2020,Otto2024}.} 

The above approaches all use linear dynamics, but there are many data-driven methods \bk{for dynamical systems} that identify non-linear dynamics.
For example, one may identify governing equations \bk{as} a sparse subset of candidate non-linear functions (Sparse Identification of Non-linear Dynamics, or \emph{SINDy})~\cite{Brunton2016sindy}.
Neural ODEs~\cite{Chen2018} use an ODE $\tfrac{\dd \x}{\dd t} = f(\x,t;\param)$ where the function $f$ is learnt from training data.
Neural operators directly learn mappings between function spaces~\cite{Lu2021,Kovachki2023,Li2020fno} and neural \bk{processes~\cite{Garnelo2018}} define probabilistic encoder--decoder \bk{architectures, and can be} used for modelling dynamical \bk{systems}\bk{~\cite{Alkin2024}}.
Even large language models have been explored for modelling and prediction of dynamical systems~\cite{liu2024}.

The {probabilistic} foundation of our work is the use of \bk{a}  proper scoring rule for optimisation \bk{to provide an inexpensive loss function targeting sharp, well-calibrated forecasts}. %
Proper scoring rules assign a score to a probability distribution and \bk{sample} pair.
The expected value of the score from a proper scoring rule is maximised when the probability distribution matches the target \bk{data} distribution~\cite{Gneiting2007}.
\bk{Therefore,} optimising a model using the expected score from a proper scoring rule targets a probabilistic model that is well calibrated with respect to the target data distribution.
We will use the Continuous Ranked Probability Score \bk{(CRPS)~\cite{Gneiting2007}.}
\bmk{Training with CRPS is the state of the art in probabilistic modelling, providing more accurate and less costly models than, e.g., variational auto-encoders, flow matching and diffusion~\cite{Ho2020,Price2025}, in applications in PDE simulation~\cite{diaconu2026probabilisticretrofittinglearnedsimulators} and weather forecasting~\cite{Alet2025,Scholz2026}.}
Notably, despite the popularity of CRPS, it is defined for scalar random variables and must be aggregated for multivariate state estimates, ignoring covariance structure.
\bk{Aggregated CRPS has still proven to be effective in these works, and \bmk{the score} is cheap to \bmk{evaluate at training time}, requiring only a small number of model samples.}
To include the covariance structure, kernel-based scoring rules~\cite{Gneiting2007,Bouchacourt2016} and test aggregation~\cite{Schrab2023} have been explored in the literature but are beyond the scope of this work.
\bmk{Deep ensembles are a popular alternative for uncertainty quantification in deep learning~\cite{Lakshminarayanan2017}.
These provide a powerful but expensive approximation that requires training several independent models.
A deep ensemble will later be used as a baseline in our simpler numerical experiments.
}

\bk{Our approach learns a probabilistic forecasting model built from Koopman motivated observation--linear dynamics--state reconstruction, and is trained to minimise CRPS.}
Existing DMD and Koopman-type methods generally minimise mean square approximation error with respect to the data.
A small number of probabilistic methods~\cite{Morton2019,Pan2020,Otto2024} use variational inference \bk{or expectation--maximisation} to fit the unknown data distribution.
\bmk{Such variational models are known to be challenging to fit successfully~\cite{Turner2011}.}
To our knowledge, no existing Koopman-type methods \bk{use sample-based models that are} optimised via the expected score from a proper scoring rule, targeting well-calibrated and interpretable forecasting.

The remainder of the paper is structured as follows: \Cref{sec:background} provides an overview of Koopman-type approximation techniques, and probabilistic modelling with a focus on proper scoring rules.
\bk{Then, the contributions are:}
\begin{enumerate}
    \item In \Cref{sec:method}, a novel probabilistic Koopman-mode approximation architecture is outlined.
    This uses conventional deep learning architectures such as multi-layer perceptrons and convolutional neural networks that define implicit distributions via the pushforward of Gaussian noise samples.
    The probabilistic observation--reconstruction functions are paired with a structured linear dynamics \bk{operator} to target identification of eigenfunctions of the Koopman operator.
    \item In \Cref{sec:method.training}, a \bk{novel CRPS-based loss and} curriculum-based training scheme is presented.
    The resulting model can be sampled to produce future state predictions capturing the predictive uncertainty arising from noisy data, and the variability in the identified dynamics.
    \item In \Cref{sec:results}, extensive numerical experiments demonstrate that the learnt observable functions can indeed approximate known Koopman eigenfunctions, that multi-modal forecast distributions can be produced, and that the architecture extends to practical fluid dynamics test problems with high dimensional state spaces.
    \item \bk{Comparison with a deterministic formulation shows the strengths of a calibrated probabilistic model, and evaluation against a deep ensemble~\cite{Lakshminarayanan2017} shows our model's comparable quality and greater expressivity at reduced computational cost.}
\end{enumerate}

\section{Background}\label{sec:background}

The linearisation of non-linear dynamical systems is first introduced via the Koopman operator, followed by an outline of Dynamic Mode Decomposition and its extensions.
\bk{Then,} probabilistic forecasting is described with a focus on proper scoring rules, and in particular the definition of the Continuous Ranked Probability Score (CRPS) as an optimisation loss.
In \Cref{sec:method}, the combination of the Koopman-type \bk{methods} and probabilistic forecasting with proper scoring rules will produce a probabilistic approximation method that provides both sharp forecasts and uncertainty quantification.

\subsection{Linear Dynamics for Non-Linear Dynamical Systems}\label{sec:background.koopman}

Consider a discrete-time autonomous dynamical system $(\X,\F)$ where $\X$ is the state space, that evolves through the (potentially non-linear) dynamics $\F:\X \to \X$,
\begin{equation}\label{eq:background.system}
    \x_{k+1} = \F(\x_k) \text{ where time index $k\in \N_0 := \N \cup \{0\}$ and $\x_0\in\X$.}
\end{equation}
Although most physical systems are continuous in time, observation and computation occur in discrete timesteps, so \eqref{eq:background.system} is a natural setting for our data-driven approximation.
The state space $\X$ may be infinite-dimensional, and in this case we consider finite-dimensional states arising from a numerical discretisation.
\bk{For PDEs, like those considered later in the numerical experiments, the state may be a function defined on a two-dimensional spatial domain.
A discretised state $\x$ may consist of nodal values on a rectangular spatial mesh with $\nX = \nX_1 \times \nX_2$ nodes,
and the discretised state may be viewed as an element of either $\R^{\nX}$ or $\R^{\nX_1 \times \nX_2}$, as appropriate.}

To identify linearised dynamics for a general non-linear problem of the form \eqref{eq:background.system}, the Koopman operator is \bk{introduced~\cite{Koopman1931,Koopman1932}.}
It is the operator that describes the evolution of \emph{observables} $g:\X\to\C$ \bk{acting on the state} of a dynamical system \eqref{eq:background.system}. %
\bk{Here and throughout, `observation' refers to this lifting of the state $\x\in\X$ by observables $g$ to latent values $g(\x)\in\C$, rather than to recorded state data.}
For an appropriate space of observables $\Gspace$, \bk{often the \bk{space of functions} ${L^2(\X,\mu)}$ for some measure $\mu$ on the state space $\X$,} the Koopman operator for the dynamical system \eqref{eq:background.system} is the operator $\koopman:\Gspace \to \Gspace$ defined such that
\begin{equation}
    [\koopman g](\x) = g(\F(\x)) \text{ for all observables } g\in\Gspace \text{ and states }\x \in \X.
\end{equation}
See, e.g.,~\cite{Budii2012,Brunton2022,Colbrook2024} for modern expositions.
For $g_1, g_2 \in \Gspace$ and $c_1, c_2 \in \C$, it is easily shown that the Koopman operator $\koopman$ is linear provided {that} $\Gspace$ is a vector space,
    $\koopman(c_1 g_1+c_2 g_2)(\x) = c_1 [\koopman g_1](\x) + c_2 [\koopman g_2](\x)$.

Our method will seek approximations to eigenfunctions of the Koopman operator $\koopman$ and build a modal approximation \bmk{to the dynamics}.
That is, for a function space $\Gspace$ with appropriate norm $\Vert \cdot \Vert_{\Gspace}$, the method will try to identify functions $\eigfn\in\Gspace$ such that,
\begin{equation}\label{eq:background.eigenfunction}
    \koopman \eigfn = \discreteeig \eigfn, \quad \Vert \eigfn \Vert_{\Gspace} = 1, \quad \lambda \in \C.
\end{equation}
For example, an appropriately normalised constant function $\eigfn(\x) \equiv C$ for all $\x\in\X$ satisfies \eqref{eq:background.eigenfunction} with $\lambda=1$ --- however $\eigfn$ tells nothing of the system state.
The existence of additional eigenfunctions is not guaranteed and the dynamics may be associated with a continuous spectrum~\cite{Mezic2005,colbrook2024limitspowerskoopmanlearning}.
However, for the purpose of this work, the following working assumption is made: the dynamics of \eqref{eq:background.system} can be accurately captured in a finite-dimensional subspace defined by the span of eigenfunctions of the Koopman operator, as defined in \eqref{eq:background.eigenfunction}. 

{Another} issue is that a data-driven method may only identify $\epsilon$-pseudo-eigenfunctions \bk{rather than true eigenfunctions.
That is, a data-driven algorithm may identify $g\in\Gspace$ satisfying,}
\begin{equation}\label{eq:koopman.pseudoeigenfunction}
    \Vert \koopman g - \lambda_\epsilon g \Vert_{\Gspace} \leq \epsilon,  \quad \Vert g \Vert_{\Gspace} = 1, \quad \lambda_\epsilon \in \C,
\end{equation} 
instead of \bk{the exact eigenvalue equation \eqref{eq:background.eigenfunction}.}
Identified pseudo-eigenfunctions may introduce spurious eigenvalues \bk{which can} corrupt an \bk{approximation~\cite{Colbrook2023,Colbrook2023res}.}
For our purpose, it is assumed that the span of identified functions define{s} an appropriate subspace to capture the system dynamics, and moving forward we will refer to identified approximate eigenfunctions as eigenfunctions.

After identifying a set $\{\eigfn_i\}_{i=1}^{\nZ}$ of $\nZ$ linearly independent eigenfunctions, for any function $g$ in the space $\Gspace_{\eigfn}:=\spanop\{\eigfn_i\}_{i=1}^{\nZ}$ {then}
    $g(\x):= \g^\top \Eigfn(\x) \text{ for some } \g \in \C^{\nZ}$,
where $\Eigfn(\x):=[\eigfn_1(\x), \eigfn_2(\x), \dots, \eigfn_\nZ(\x)]^\top$.
Defining $\K := \diagop(\lambda_1,\lambda_2,\dots,\lambda_\nZ) \in\C^{\nZ \times \nZ}$, then
\begin{equation} \label{eq:kmd}
    g(\F(\x))=[\koopman g](\x) = \g^\top \bk{\K} \Eigfn(\x).
\end{equation}
\bmk{Equation~\eqref{eq:kmd}} gives a linear representation $\K$ of the dynamics in the space $\Gspace_{\eigfn}$, with non-linear observation $\Eigfn:\X \to \C^\nZ$ and linear reconstruction of a quantity of interest $g$.

In practice, many simultaneous quantities of interest {may be gathered} into a vector-valued function $\vec{g}:\X \to \C^{\nX}$.
If $\vec{g}:= \G \Eigfn$, then the Koopman mode decomposition (KMD) $(\Eigfn,\K,\G)$ is
\begin{equation} \label{eq:background.kmd}
    \vec{g}(\F^k(\x)) = \G \K^k \Eigfn(\x),
\end{equation}
where the columns $\g_i \in \C^{\nX}$ of $\G=[\g_1, \g_2, \dots, \g_\nZ]\in\C^{\nX\times\nZ}$ are the Koopman modes, representing coherent structures in vector-valued observables that are associated to a particular eigenfunction $\eigfn_i$.
The KMD \eqref{eq:background.kmd} is the starting point for our approximation: $\Eigfn:\X \to \C^\nZ$ is the state observation, the diagonal matrix $\K$ defines the linear dynamics, and the linear operator $\G$ is the reconstruction operator for a vector valued function $\vec{g}\in \left(\Gspace_{\Eigfn}\right)^{\nX}$.

In practice, data-driven methods are used to identify observation functions, linear dynamics and reconstruction operators.
\bk{An overview of existing methods is listed in \Cref{table:background.summary} and the remainder of this subsection will outline the key features of each class of methods.
\begin{table}\centering
\footnotesize
    \caption{Selected development timeline for DMD, and Koopman-type approximation algorithms.}
    \label{table:background.summary}
    \begin{tabular}{l l l l l l}
        \toprule
        \textbf{Reference}              &
        \textbf{Year} &
        \textbf{Observables} &
        \textbf{Dynamics}   &
        \textbf{Reconstruction} &
        \textbf{Data fit} \\
        \midrule
        \emph{DMD}&&&&&\\
        Rowley et al.~\cite{Rowley2009} & 2009 & Linear  & Linear  & Linear  & MSE \\
        Schmid~\cite{Schmid2010} & 2010 & Linear  & Linear  & Linear  & MSE \\
        \midrule
        \emph{Extended DMD}&&&&&\\
        Williams et al.~\cite{Williams2015,Williams2015b} & 2015 & Non-linear & Linear  & Linear  & MSE  \\
        Li et al.~\cite{Li2017} & 2017 & NN & Linear & Linear  & MSE  \\
        \midrule
        \emph{Machine Learning}&&&&&\\
        Takeishi et al.~\cite{Takeishi2017b} & 2017 & NN & Linear & NN  & MSE \\
        Lusch et al.~\cite{Lusch2018} & 2018 & NN & Diagonal & NN  & MSE\\
        Otto et al.~\cite{Otto2019} & 2019 & NN & Diagonal & NN  & MSE\\
        Erichson et al.~\cite{Erichson2019}    & 2019 & NN     & Linear NN     & NN     & MSE   \\
        Guo et al.~\cite{Guo2025} & 2025 & NN & Parametric & Linear  & MSE \\
        \midrule
        \emph{Probabilistic}&&&&&\\
        Takeishi et al.~\cite{Takeishi2017} & 2017 & Linear Gaussian  & Diagonal  & Linear Gaussian & MCMC \\
        Morton et al.~\cite{Morton2019} & 2019 & Gaussian (NN) & Linear & Gaussian (NN) & VI \\
        Pan et al.~\cite{Pan2020} & 2020 & Bayesian NN & Constrained & Bayesian NN & VI \\
        Otto et al.~\cite{Otto2024} & 2024 & -- & Linear Gaussian & Gaussian (NN)  & E--M \\
        \textbf{This work} & -- & NN (probabilistic) & Diagonal & NN (probabilistic) & CRPS\\
        \bottomrule
    \end{tabular}
\end{table}}
To fix notation, the training dataset $\data$ consists of $\nData$ trajectories, %
\begin{equation}\label{eq:background.data}
\data = \left\{\data^{(i)}\right\}_{i=1}^\nData =\left\{\left(\gt{\x}_0^{(i)}, \gt{\x}_1^{(i)}, \dots, \gt{\x}_T^{(i)}\right)\right\}_{i=1}^{\nData},
\end{equation}
where $\gt{\x}_k^{(i)}=F^k(\gt{\x}_0^{(i)}) \in \X$ is the ground truth state at time $k$ for trajectory $i$ with initial state $\gt{\x}_0^{(i)}$.

\emph{Dynamic mode decomposition} (DMD) provides a tractable method for computing an approximation of the same structure as KMD \eqref{eq:background.kmd} for state spaces $\X \subseteq \C^\nX$~\cite{Rowley2009,Tu2014,Schmid2010}.
Observation and reconstruction are based upon linear projection $P^\top:\X\to\C^{\nZ}$, and a linear operator $\K \in L(\C^{\nZ} ; \C^{\nZ})$ is found by minimising the least-squares approximation error between consecutive state data,
\begin{equation}
\label{eq:background.dmd.simple}
    \K_{DMD} = \argmin_{\K \in L(\C^{\nZ};\C^{\nZ})} \sum_{i=1}^{\nData} \Vert \gt\x_1^{(i)} - P \K P^{\top} \gt\x_0^{(i)} \Vert^2,
\end{equation}
where trajectories of length $T>1$ can be split into consecutive pairs to form the complete dataset.
Diagonalisation of the dynamics results in a representation that shares the characteristics of the Koopman mode decomposition: \bk{system states $\x$ are observed via linear projection onto the DMD modes, dynamics are captured via the eigenvalues of the matrix $\K_{DMD}$, and the state vector $\x$ is reconstructed as a linear combination of the DMD modes.}
Additional structure can be enforced, for example, by constraining the manifold of matrices in the least-squares approximation \eqref{eq:background.dmd.simple}, parametric and control problems can be analysed by concatenating parameters \bk{to the state}~\cite{Proctor2016,Huhn2023,Andreuzzi2023}, and uncertainty estimates can be produced using an ensemble of approximations~\cite{Sashidhar2022}.

Generally, the mean square error \eqref{eq:background.dmd.simple} \bk{cannot} be made zero as the \bk{underlying dynamics are non-linear}.
Instead, one may choose to use \emph{Extended DMD} (EDMD)~\cite{Williams2015} and seek linear dynamics for a dictionary of linearly independent observables $\{\dict_i\}_{i=1}^{\nZ}$, $\dict_i:\X \to \C$ for $i\in\{1,2,\dots,\nZ\}$.
Denote the space spanned by the observables $\Gspace_{\dict}:=\spanop\{\dict_i\}_{i=1}^{\nZ} \subset \Gspace$.
The goal is to identify linear dynamics $\K$ minimising the residual, 
\begin{equation}\label{eq:background.edmd.res}
    {\residual(K) := \sum_{j=1}^\nZ} \residual_j(K) :=  {\sum_{j=1}^\nZ} \left\Vert \dict_j(\F(\cdot)) - \sum_{i=1}^{\nZ} \dict_i \K_{i,j} \right\Vert_{\Gspace}
\end{equation}
for {an} appropriate norm on $\Gspace$ \bk{\cite{Williams2015,Colbrook2025}.
EDMD uses the right multiplication convention for dynamics.}
In the data-driven setting, \eqref{eq:background.edmd.res} is minimised by treating the training data as a quadrature method over the domain $\X$.
Omitting the details~\cite{Williams2015,Colbrook2024,Colbrook2025}, the resulting method minimises \eqref{eq:background.edmd.res} through the dynamics $\K$ computed via matrices $\Psi_{X},\Psi_{Y} \in \C^{\nData \times \nZ}$,
\begin{subequations}\label{eq:background.edmd.K}
    \begin{align}
    [\Psi_X]_{i,j}&:=\dict_j\left(\gt\x_0^{(i)}\right), \quad[\Psi_Y]_{i,j}:=\dict_j\left(\gt\x_1^{(i)}\right),\\ \K := \psinv{\Gram}A &:= \psinv{(\Psi_X^* W \Psi_X)} \Psi_X^* W \Psi_Y,
    \end{align}
\end{subequations}
where $W=\diagop(w_1,w_2,\dots,w_\nData)$ is a diagonal matrix of quadrature-like weights.
For a target function $g\in\Gspace_{\dict}$, then $g(\x):=\Dict^{\top}(\x)\g$ for some $\g\in\C^{\nZ}$ where the observables are collated as $\Dict^{\top}(\x)=[\dict_1(\x),\dict_2(\x),\dots,\dict_\nZ(\x)]$ and $g$ is then approximated by $\estimate{g}_k : = \Dict(\gt\x_0)^{\top}\K^{k}\g \approx g(\gt\x_k)$.
Diagonalisation of the dynamics gives a representation with the characteristic structure of the Koopman mode decomposition \eqref{eq:background.kmd}: \emph{non-linear} observation, diagonalised linear dynamics, and linear observable reconstruction.

\bk{The choice} of a suitable dictionary is essential for EDMD, but \bk{selecting this is difficult,} particularly for high (or infinite) dimensional state spaces $\X$ such as those arising from PDE discretisation.
Instead, problem specific observables can be handcrafted but a more general solution is to learn them.
For example, the observables $\Dict:\X\to\C^{\nZ}$ can be parameterised neural \bk{networks} $\Dict(\x):=\Dict_{NN,\param_{\dict}}(\x)$ that are optimised for fixed matrix $\K$ under a mean square error loss, with $\K$ \bk{periodically re}computed via the EDMD method \eqref{eq:background.edmd.K}~\cite{Li2017}.
Variations on this approach are taken in~\cite{Takeishi2017b,Yeung2019,Otto2019,Erichson2019,Lusch2018}, including directly targeting diagonalisation of the dynamics $\K$ to learn a dictionary of eigenfunctions and eigenvalues of the Koopman operator~\cite{Lusch2018}.

\bk{Predicted quantities of interest often include the state itself.
For linear reconstruction, this requires the identity functions $g_{x_i}(\x) = x_i$ for each state variable to belong to the span of the EDMD observables $\Gspace_{\dict}$.
\bmk{Including these functions} \bk{can be} restrictive as \bk{it forces} linear dynamics to be identified for the state variables, which may be hard, or even impossible, with a finite dictionary~\cite{Brunton2022}.
More recent Koopman-type approaches \bk{instead} use a trainable non-linear neural network {reconstruction} $g_{\param_{g}}:\Z\to\X$, $g_{\param_{g}}(\z):=g_{NN,\param_{g}}(\z)$ with trainable parameters $\param_{g}$ optimised with respect to a least-squares loss~\cite{Takeishi2017b,Yeung2019,Otto2019,Erichson2019,Lusch2018}.}

In~\cite{Lusch2018,Otto2019} the dynamics operator is constrained using block rotation matrices, providing a real-valued representation of the complex-valued observation space dynamics.
Observation pairs seek approximations to the Koopman eigenfunctions $\eigfn$ defined in \eqref{eq:background.eigenfunction}.
That is, a set of $\nZ\in2\N$ observables are paired as $\{(\dict_{2i-1},\dict_{2i})\}_{i=1}^{\nZ/2}$ to satisfy,
\begin{subequations}
\begin{align}
    \left[\dict_{2i-1} + \im \dict_{2i}\right](\F(\x)) &= \koopman \left[\dict_{2i-1} + \im \dict_{2i}\right](\x) = \lambda_i \left(\dict_{2i-1} + \im \dict_{2i}\right)(\x), \label{eq:background.realrepresentation}\\
    \left[\dict_{2i-1} - \im \dict_{2i}\right](\F(\x)) &= \koopman \left[\dict_{2i-1} - \im \dict_{2i}\right](\x) = \conjugate{\lambda}_i \left(\dict_{2i-1} - \im \dict_{2i}\right)(\x). \label{eq:background.conjugate}
\end{align}
\end{subequations}
where for real-valued data, the eigenfunctions and eigenvalues come as complex conjugate pairs.
With $\lambda_i = \exp(\mu_i \Delta t)\exp(\im \theta_i \Delta t)$ for a fixed $\Delta t \in \R$,
\eqref{eq:background.realrepresentation} and \eqref{eq:background.conjugate} are summarised by,
\begin{subequations}\label{eq:background.koopman.rotation}
\begin{align}
    R(\theta_i \Delta t) &:= \begin{bmatrix} \cos(\theta_i \Delta t)& -\sin(\theta_i \Delta t)\\ \sin(\theta_i \Delta t)& \cos(\theta_i \Delta t) \end{bmatrix}, \label{eq:background.koopman.rotation.R}\\
    \begin{bmatrix} \dict_{2i-1}(\x_{k+1}) \\ \dict_{2i}(\x_{k+1}) \end{bmatrix} &\approx e^{\mu_i \Delta t} R(\theta_i \Delta t) \begin{bmatrix} \dict_{2i-1}(\x_k) \\ \dict_{2i}(\x_k) \end{bmatrix}. \label{eq:background.koopman.rotation.approx}
\end{align}
\end{subequations}
The rotation blocks are combined to give block diagonal dynamics,
\begin{equation}\label{eq:background.K.rotation}
\K := {\small\begin{bmatrix}
    e^{\mu_1 \Delta t} R(\theta_1 \Delta t) & \mathbf{0}                             & \cdots & \mathbf{0} \\
    \mathbf{0}                              & e^{\mu_2 \Delta t} R(\theta_2 \Delta t) & \cdots & \mathbf{0} \\
    \vdots                                  & \vdots                                  & \ddots & \vdots     \\
    \mathbf{0}                              & \mathbf{0}                              & \cdots & e^{\mu_{\nZ/2} \Delta t} R(\theta_{\nZ/2} \Delta t)
\end{bmatrix}} \in \R^{\nZ \times \nZ},
\end{equation}
resulting in the linearised observable dynamics $\Dict(\x_{k+1}) = \K\,\Dict(\x_k)$ {for the full dictionary $\Dict$}.
In~\cite{Lusch2018}, \bmk{continuous spectra are targeted with state-dependent} eigenvalues $\{\lambda_i(\z)\}_{i=1}^{\nZ/2}$ that are optimised along {with} the dictionary $\Dict$, with $\Lambda(\z):=\{\lambda_i(\z)\}_{i=1}^{\nZ/2} = \Lambda_{NN,\param_\K}(\z)$ defining {the} dynamics matrix $\K(\z):=\K(\Lambda(\z))$.
The real-valued formulation \eqref{eq:background.koopman.rotation} will be the starting point of our method in \Cref{sec:method}, with a trainable approximation structured as
\begin{equation}\label{eq:background.kmd.ml}
    \estimate{\x}_k := g_{\param_g}(\K^k{(\Lambda_{NN,\param_\K})} \Dict_{\param_\Dict}(\x_0)) \approx \F^k(\x_0).
\end{equation}

Further to this, a probabilistic formulation will be used for each of the three components.
Bayesian DMD~\cite{Takeishi2017}, which builds on probabilistic PCA~\cite{Tipping1999}, places priors on DMD modes and eigenvalues and quantifies the approximation uncertainty via Gibbs sampling.
In Deep Variational Koopman approximation~\cite{Morton2019}, the observables are latent random variables with parameterised Gaussian distributions inferred by maximising the evidence lower bound (ELBO).
Sampling of the observables yields an ensemble of linear dynamics models and hence uncertainty-quantified state predictions.
In~\cite{Pan2020} the observables are learnt as neural network non-linear corrections to the DMD projection.
A Gaussian mean-field variational posterior is inferred over the neural network weights and a structured tri-diagonal dynamics matrix.
Predictions are made by sampling the posterior.
In~\cite{Otto2024}, a Koopman-type approximation is fit via a linear-Gaussian model and an expectation--maximisation algorithm.
These methods are fit via likelihood-based \bk{inference.}
In contrast, our approach in \Cref{sec:method} will optimise an implicit distribution via a proper scoring rule, which will be outlined in \Cref{sec:background.scoring}.

\subsection{Probabilistic Forecasting and CRPS}
\label{sec:background.scoring}
We outline probabilistic forecasting via proper scoring rules, with a focus on the Continuous Ranked Probability Score (CRPS).
\bmk{CRPS} will be used \bmk{in} a loss function in \Cref{sec:method} to optimise a probabilistic head defined via the pushforward of a low dimensional Gaussian noise sample.

The goal of this work is to provide a sharp but well-calibrated probabilistic state forecast that is consistent with the training data statistics of a future event~\cite{Gneiting2007}.
The quality of a probabilistic forecast can be measured by a \emph{scoring rule}.

\begin{definition}[Strictly proper scoring rules~\cite{Gneiting2007}]
\label{def:background.scoring.proper}
Consider the sample space $\Omega$ and Borel probability measures $\mathcal{P}$ on $\Omega$.
A scoring rule $s:\mathcal{P}\times \Omega \to \R$ is a function mapping a probability measure $P \in \mathcal{P}$ and event $\omega \in \Omega$ to a \emph{score}.
The expected score $S(\cdot, Q):\mathcal{P} \to \R$ for a reference measure $Q$ is,
\begin{equation}
    S(P,Q):=\E_Q[s(P,\cdot)]=\int_{\Omega} s(P,\omega) dQ(\omega).
\end{equation}
A scoring rule $s$ is proper relative to $\mathcal{P}$ if for all $P,Q\in\mathcal{P}$ the expected score $S$ satisfies $S(Q,Q) \geq S(P,Q)$.
A scoring rule $s$ is strictly proper if it is proper, and the expected score  $S(Q,Q)=S(P,Q)$ if and only if $P = Q$.
A more precise measure-theoretic definition can be found in~\cite{Gneiting2007}.
\end{definition}

A popular choice of proper scoring rule is the \emph{Continuous Ranked Probability Score} (CRPS) on sample spaces $\Omega \subseteq\R$~\cite{Hersbach2000,Gneiting2007}.
\bmk{CRPS} is normally written as a non-negative function to be minimised, and this is the convention we will use \bk{throughout}.
\begin{definition}[Continuous Ranked Probability Score]
\label{def:background.crps}
For a probability measure $P$, the CRPS is defined using the cumulative distribution function $F_P$ and a \bk{target datum} $\omega$ by,
\begin{equation}\label{eq:crps.definition}
    \crps(F_P,\omega)=  \int_{\R} \left(F_P(y) - \indicator_{y\geq\omega}\right)^2 \dd y,
\end{equation}
where $\indicator_{y\geq\omega}$ is the indicator function that is 1 if $y \geq \omega$ and 0 otherwise.
\end{definition}
CRPS defines a strictly proper score (when negated) with respect to probability measures with finite first moment~\cite{Gneiting2007}, and in this case, \eqref{eq:crps.definition} can be alternatively expressed using random variables $X,X' \sim P$ as
\begin{equation}\label{eq:background.crps.alt}
    \crps(F_P, \omega) = \E_P[|X - \omega|] - \frac{1}{2} \E_P[|X- X'|].
\end{equation}
\bk{In the formulation \eqref{eq:background.crps.alt}, two independent samples $X$ and $X'$ are drawn from the forecast distribution $P$.
The first term, $\E_P[|X - \omega|]$, penalises the absolute error.
The second term, $-\frac{1}{2} \E_P[\lvert X - X' \rvert ]$, rewards ensemble spread.
CRPS will penalise an over-dispersed distribution $P$ via the first term, whilst an over-confident approximation will receive little reward from the second term.}
If the forecast distribution $P$ is a point mass, then \eqref{eq:background.crps.alt} reduces to the absolute error.
\bk{In fact, the formulation~\eqref{eq:background.crps.alt} is a specific example of a kernel score~\cite{Gneiting2007} with kernel $\kernel(u,v)=\lvert u - v \rvert$,
\begin{equation}\label{eq:background.kernelscore}
    \score(P, \omega) = \E_P[\kernel(X,\omega)] - \frac{1}{2} \E_P[\kernel(X,X')].
\end{equation}
More generally, $\kernel:\Omega\times\Omega\to\R$ can be any \bmk{conditionally} negative definite kernel.
}
An unbiased estimator of~\eqref{eq:background.crps.alt} is~\cite{Zamo2018},
\begin{equation}\label{eq:background.crps.estimator}
    \estimate{\crps}(\samples,\omega) := \frac{1}{\nSamples} \sum_{i=1}^{\nSamples} \lvert x^{(i)} - \omega \rvert - \frac{1}{2 \nSamples(\nSamples-1)} \sum_{i=1}^{\nSamples}\sum_{j=1}^{\nSamples} \lvert x^{(i)} - x^{(j)} \rvert,
\end{equation}
where $\samples:=\{x^{(i)}\}_{i=1}^{\nSamples}$ is a \bk{set of independent samples $x^{(i)} \sim P$}.

Minimising CRPS is a popular optimisation strategy in probabilistic forecasting, e.g.{,} in applications in weather~\cite{Alet2025,Lang2026,Scholz2026} and PDEs~\cite{diaconu2026probabilisticretrofittinglearnedsimulators}.
\bk{It is simple to interpret, and the estimator \eqref{eq:background.crps.estimator} is cheap to compute, with practitioners often using as few as $\nSamples=2$ samples per data point during training~\cite{Alet2025}.}

In practice, \bk{datasets are} often multivariate.
To this end, \eqref{eq:crps.definition} is applied to each \bk{state component} independently, and the resulting scores are averaged across \bk{components.}
That is, for random vectors $X,X' \sim P$ with $\nX$ \bk{components,}
\begin{equation}\label{eq:background.crps.multivariate}
    \estimate{\crps}(\samples, \vec{\omega}) = \frac{1}{\nX} \sum_{j=1}^{\nX} \estimate{\crps}(\samples_j, \omega_j),
\end{equation}
for model samples $\samples_j$ and target data $\omega_j$ for the $j$-th \bk{component.}
\bmk{Equation~\eqref{eq:background.crps.multivariate}} does not capture dependencies between \bk{components,} but is effective in practice~\cite{Alet2025,diaconu2026probabilisticretrofittinglearnedsimulators}.
Multi-variate scoring \bk{rules} such as \bk{the energy and other kernel scores} are not considered in this work~\cite{Gneiting2007,Bouchacourt2016,Schrab2023}.

\section{Method}
\label{sec:method}
Our probabilistic approximation strategy is built upon the structure of the Koopman mode decomposition defined in \eqref{eq:background.kmd} \bk{and CRPS based training.}
\bk{First, the probabilistic model is defined in \Cref{sec:method.architecture}, followed by the training strategy in \Cref{sec:method.training}.}
\subsection{\bk{Probabilistic model}}
\label{sec:method.architecture}
\begin{figure}\centering
    \includegraphics{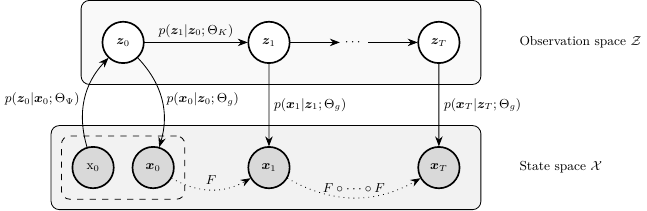}
    \caption{Graphical model for the probabilistic KMD forecasting method.}
    \label{fig:method.graphicalmodel}
\end{figure}

The three components are:
\begin{itemize}
    \item a probabilistic observation strategy $\sample{\z} \sim \Dict_{\param_\Dict}(\x)$, that will be optimised to identify pairs of observables approximating a real-valued representation of complex-valued eigenfunctions $\eigfn$, i.e., satisfying \eqref{eq:background.realrepresentation},
    \item a random linear dynamics operator $\K \sim \K_{\param_\K}$, $\K: \Z \to \Z$ with the block diagonal structure \eqref{eq:background.K.rotation} corresponding to the eigenvalues \bk{$\{\lambda_i\}_{i=1}^{\nZ/2}$} associated with the identified observables $\Dict_{\param_\Dict}$,
    \item and a probabilistic \emph{non-linear} reconstruction strategy $\sample{\x} \sim g_{\param_g}(\z)$ that will be optimised to reconstruct a state $\x \in \X$ conditioned on an observation $\z\in\Z$.
\end{itemize}
Given an initial state $\gt{\x}_0 \in \X$, the observation--dynamics--reconstruction model is written via conditional distributions,
\begin{equation}\label{eq:sampling_conditionals}
    \begin{aligned}
    \text{observation } \quad &p(\z_0 | \gt{\x}_0; \param_\Dict),& \\
    \text{linear dynamics }\quad &p(\K ; {\param_\K}),& \\
    \text{observable propagation }\quad &p(\z_k | \z_{k-1},\K) = \delta(\z_k - \K \z_{k-1}), &  \quad k = 1, \ldots, T, \\
    \text{and reconstruction }\quad &p(\x_k | \z_k; \param_g), &\quad k = 0, \ldots, T,
\end{aligned}
\end{equation}
with joint predictive distribution
\begin{equation}\label{eq:method.probmodel}
    \begin{aligned}
    p(\x_{0:T} | \gt{\x}_0; \param) = \int_{\Z^{T+1}} \int_{\R^{\nZ \times \nZ}} & p(\z_0 | \gt{\x}_0; \param_\Dict) \, p(\K ; \param_\K) \\ & \times \prod_{k=0}^T p(\x_k | \z_k; \param_g) \prod_{k=1}^T p(\z_k | \z_{k-1}, \K) \, \dd \K \dd \z_{0:T}.
    \end{aligned}
\end{equation}
Samples are generated using a \bk{discrete-time} rollout \bk{$\z_k:= \K^k\z_0$} in the observation space $\Z$ with a fixed dynamics operator $\K$ per {trajectory} sample.
\bmk{The initial state} $\x_0$ is included in the predictive distribution as we consider the state reconstruction from initial \bmk{data} $\gt{\x}_0$ {after observation}.
\bmk{Equation~\eqref{eq:method.probmodel}} is illustrated in the graphical model in \Cref{fig:method.graphicalmodel}.
The model parameters $\param = [\param_\Dict, \param_\K, \param_g]$ will be trained by optimising with respect to a proper scoring rule based loss, namely CRPS as described in \eqref{eq:background.crps.multivariate}, \bk{using} training data $\data$.

\begin{figure}\centering
\includegraphics{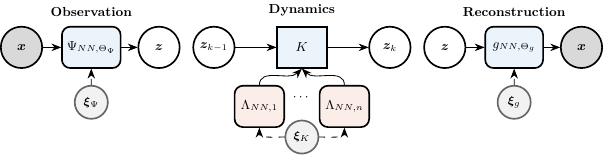}
    \caption{Diagram of the model blocks (observation, linear dynamics, and reconstruction) with noise injection mechanism. Each block uses an independent noise sample, respectively, $\noise_\Dict\sim\Ndist(\vec{0},I_{\nNoise_\Dict})$, $\noise_\K\sim\Ndist(\vec{0},I_{\nNoise_\K})$, $\noise_g\sim\Ndist(\vec{0},I_{\nNoise_g})$.}
    \label{fig:noise_injection}
\end{figure}

The probabilistic observation function $\z\sim\Dict_{\param_\Dict}(\x)$ is defined by a neural network $\Dict_{NN,\param_\Dict}:\X \times \R^{\nNoise_\Dict} \to \Z$ that maps a state $\x\in\X$ paired with a sample of Gaussian noise $\noise_\Dict \sim \Ndist(\vec{0}, I_{\nNoise_\Dict})$, where $\nNoise_\Dict\in\N$, to an observation $\z \in \Z$.
That is, the conditional distribution $p(\z \given \x; \param_\Dict)$ is the pushforward of $\Ndist(\vec{0},I_{\nNoise_\Dict})$ under the map $\Dict_{NN}(\x,\cdot)$.
This approach is common and successful in probabilistic forecasting~\cite{Alet2025,diaconu2026probabilisticretrofittinglearnedsimulators}.
In the numerical results in \Cref{sec:results}, this will \bmk{be implemented as} either noise concatenated to the state vector as input into a multi-layer perceptron, or noise Feature-wise Linear Modulation (FiLM) conditioning~\cite{Perez2018} in a convolutional neural network for structured state spaces, e.g., when representing PDE solutions.
The specific architecture used in the numerical experiments is described in \Cref{sec:architecture} whilst a sketch of the architecture is provided in \Cref{fig:noise_injection}.
\bmk{The pushforward construction} is simpler than, e.g., a mean field Bayesian neural network~\cite{Pan2020} which requires a noise sample for each network weight.

The dynamics use the block diagonal structure \eqref{eq:background.K.rotation} to represent complex-valued eigenpairs.
The eigenvalues are computed from independently parameterised neural networks that share a common input noise sample \bk{$\noise_\K\sim\Ndist(\vec{0},I_{\nNoise_\K})$.}
That is, for each eigenvalue, a scaling and rotation are computed via a neural network $\Lambda_{NN,i}:\R^{\nNoise_\K} \to \R^2, \noise_\K\mapsto(\mu_i, \theta_i)$ representing each discrete-time eigenvalue $\lambda_i := (1+\mu_i \Delta t) \exp(\im \theta_i  \Delta t) \approx \exp(\mu_i \Delta t) \exp(\im \theta_i  \Delta t)$ for some \bmk{fixed time scaling} constant $\Delta t$.
The eigenvalues are used to construct a sample of the structured linear operator $\K \sim\K_{\param_\K}$ via \eqref{eq:background.K.rotation}, \bmk{where $\param_\K$ collates the parameters of the eigenvalue neural networks $\{\Lambda_{NN,i}\}_{i=1}^{\nZ/2}$}.

\begin{remark}[Constraining spectral radius]
\label{rem:method.spectralradius}
    The linear dynamics $\K$ can be constrained to keep the spectral radius, that is $\rho = \max_{i \in \{1,2,\dots,\nZ/2\}} \lvert \lambda_i \rvert$, less than or equal to one.
    This has previously been considered in~\cite{Pan2020} using a structured tri-diagonal dynamics matrix, and in~\cite{Erichson2019} via an eigenvalue regularisation term.
    Control of the spectral radius $\rho \leq 1$  implies $\Vert \K \Vert_2 = \rho \leq 1$ given that the dynamics matrix \eqref{eq:background.K.rotation} is normal by construction.
    This ensures that the  magnitude of the observation $\z_k$ remains bounded by the norm of the initial observation,
    \begin{equation}
        \Vert \z_k \Vert_2 \leq \Vert \K^k \Vert_2 \Vert \z_0 \Vert_2 = \rho^k \Vert \z_0 \Vert_2 \leq \Vert \z_0 \Vert_2.
    \end{equation}
    {In practice, this is achieved by parameterising each eigenvalue $\lambda_i$ via a sampled pair $(\widetilde{\mu}_i, \theta_i)$ and using the \bk{alternative eigenvalue} function \begin{equation}\lambda_i = \exp( -\Delta t \texttt{SoftPlus}(\widetilde{\mu}_i)) \exp(\im \theta_i \Delta t),\end{equation} where $\texttt{SoftPlus}(x) = \log(1+\exp(x))$.
    This enforces a strict bound and ensures that $\lvert \lambda_i \rvert < 1$ for all $i$, so it cannot precisely capture dynamics with eigenvalues on the unit circle.
    }
{In the numerical experiments below, it was only necessary to constrain the eigenvalues for the Duffing oscillator experiments.
Elsewhere,} numerical blow up is avoided in part due to a curriculum based training strategy.
\end{remark}

\bk{The probabilistic reconstruction $g_{\param_g}$ is defined {similarly} to the probabilistic observations, i.e., via a neural network $g_{NN,\param_g}:\Z \times \R^{\nNoise_g} \to \X$ that maps an observation $\z\in\Z$ paired with a sample of Gaussian noise $\noise_g \sim \Ndist(\vec{0}, I_{\nNoise_g})$ to a state $\x\in\X$.}

\subsection{Training}
\label{sec:method.training}
\bmk{For a time horizon $T$}, we minimise the multivariate CRPS \eqref{eq:background.crps.multivariate} computed for \bk{the} \bk{state forecast}. 
\bk{This optimisation} targets the marginal distribution \bk{$p(\gt{x}_j) \approx p({x}_j;\param)$ for each state component $\gt{x}_j$} at each timestep \bk{$k=0,1,\dots,T$}.
\bmk{The learned distribution combines data and model uncertainty, and the training strategy does not distinguish between the two.}
As mentioned in \Cref{sec:background.scoring}, marginal CRPS does not directly target any covariance structure.
\bmk{However, targeting the marginals at a variety of different lead-times does constrain aspects of the joint distribution}.
\bk{Additionally, to promote consistent linear dynamics, a second CRPS-type term will be used to compare the forecast observation $p(\z_k\given \gt\x_0)$ and observation of the ground truth data $\gt\z_k \sim \Dict_{\param_\Dict}(\gt\x_k)$.}

The training dataset $\data$ consists of $\nData$ trajectories $\left\{\data^{(i)}\right\}_{i=1}^\nData$, as defined in \eqref{eq:background.data}, where $\gt{\x}_k^{(i)} \in \X$ is the ground truth state at time $k$ for trajectory $i$.
For each training trajectory $\data^{(i)}$, we construct an ensemble of $\nSamples$ approximate trajectories $\samples^{(i)}$ based upon the initial state $\gt{\x}_0^{(i)}$ and sample using the strategy described by \cref{eq:sampling_conditionals},
\begin{equation}\label{eq:method.samples.state}
\samples^{(i)} = \left\{\samples^{(i,j)}\right\}_{j=1}^{\nSamples} = \left\{\left\{\sample{\x}_0^{(i,j)},\sample{\x}_1^{(i,j)}, \dots, \sample{\x}_T^{(i,j)}\right\}\right\}_{j=1}^{\nSamples}
\end{equation}
where $\sample{\x}_k^{(i,j)}$ is the $j$-th sample for trajectory $i$ at time $k$.
Simultaneously, we obtain predicted observation samples $\samples^{\Z,(i)}$ and training data observations $\samples^{\Dict,(i)}$,
\begin{subequations}\label{eq:samples.obs.pair}
\begin{align}
\samples^{\Z,(i)} &= \left\{\samples^{\Z,(i,j)}\right\}_{j=1}^{\nSamples} = \left\{\left\{\sample{\z}_0^{(i,j)},\sample{\z}_1^{(i,j)}, \dots, \sample{\z}_T^{(i,j)}\right\}\right\}_{j=1}^{\nSamples} \label{eq:method.samples.obs}\\
\samples^{\Dict,(i)} &= \left\{\samples^{\Dict,(i,j)}\right\}_{j=1}^{\nSamples} = \left\{\left\{\gt{\z}_0^{(i,j)},\gt{\z}_1^{(i,j)}, \dots, \gt{\z}_T^{(i,j)}\right\}\right\}_{j=1}^{\nSamples} \label{eq:samples.obs}
\end{align}
\end{subequations}
\bmk{Samples} $\samples^{\Dict,(i)}$ are \bmk{drawn} from a different distribution to $\samples^{\Z,(i)}$ and only depend on the parameters $\param_{\Dict}$ \bmk{(each sample in $\samples^{\Z,(i)}$ will also depend on the dynamics parameters $\param_\K$)}.
The observation samples $\{\samples^{\Dict,(i)}\}_{i=1}^{\nData}$ are important for enforcing linear dynamics.

The loss function is based on CRPS, as described in \Cref{sec:background.scoring}.
\bk{Specifically, we use the sample-based multivariate estimator~\eqref{eq:background.crps.multivariate}.}
It is split into a state space loss $\loss_{\crps,\X}$ and a \bk{consistency term} $\loss_{\crps,\Z}$ \bmk{in observation space $\Z$}, both \bmk{computed over} a multi-step rollout.
The state space loss is defined for a time horizon $\bmk{T}$ for a trajectory $i$ by
\begin{equation}\label{eq:loss.state}\small
\loss_{\crps,\X}(\samples^{(i)}, \data^{(i)}) = \estimate\crps\!\left(\left\{\sample{\x}_0^{(i,j)}\right\}_{j=1}^{\nSamples}, \gt{\x}_0^{(i)}\right) + \sum_{k=1}^\bmk{T} \estimate\crps\!\left(\left\{\sample{\x}_k^{(i,j)}\right\}_{j=1}^{\nSamples}, \gt{\x}_k^{(i)}\right).
\end{equation}
The $k=0$ term is kept separate as it represents the observation--reconstruction pair, whilst all other terms are also influenced by the linear dynamics $\K$.
The consistency term is defined for a time horizon $\bmk{T}$ for a trajectory $i$ by
\begin{equation} \label{eq:loss.obs}\small
    \begin{aligned}
    \loss_{\crps,\Z}(\samples^{\Z,(i)}, \samples^{\Dict,(i)}) = \frac{1}{\nSamples} \sum_{\ell=1}^{\nSamples} \Bigg(& \estimate\crps\!\left(\left\{\sample{\z}_0^{(i,j)}\right\}_{j=1}^{\nSamples}, \gt{\z}_0^{(i,\ell)}\right)\\
    & + \sum_{k=1}^\bmk{T} \estimate\crps\!\left(\left\{\sample{\z}_k^{(i,j)}\right\}_{j=1}^{\nSamples}, \gt{\z}_k^{(i,\ell)}\right) \Bigg).    
    \end{aligned}
\end{equation}
\bmk{Equation~\eqref{eq:loss.obs}} must use the observation samples $\samples^{\Dict,(i)}$ defined in \eqref{eq:samples.obs} because there is no ground truth observation \bk{data}.
\bmk{The consistency term $\loss_{\crps,\Z}$ demands that the propagated and observed latent states agree, but should not be viewed as the application of a proper scoring rule as the reference distribution is not fixed and depends on $\param_\Dict$.}
Without the consistency term, the observation and reconstruction operations may influence the approximated dynamics in a non-linear way, obscuring the \bk{intended} linear observation space \bk{dynamics}.
\bk{This issue was observed in preliminary experiments.}
The \bk{consistency term} \eqref{eq:loss.obs} can also be \bk{minimised} by, for example, collapsing the observables to a constant, $\Dict_{\param_\Dict} \equiv \vec{1}$.
However, this would spoil the state reconstruction loss \eqref{eq:loss.state} as no meaningful information is propagated.
The state reconstruction loss \eqref{eq:loss.state} and the consistency \bk{term} \eqref{eq:loss.obs} are combined and \bk{averaged} over a batch of trajectories $\batch = \{ \data^{(b)} \}_{b=1}^{\nBatch}$,
\begin{equation}\label{eq:loss}
    \loss_{\X,\Z}(\param; \batch) := \bk{\frac{1}{\nBatch}\sum_{b=1}^{\nBatch}} \left(\lossweight_\X \loss_{\crps,\X}(\samples^{(b)}, \data^{(b)}) + \lossweight_\Z \loss_{\crps,\Z}(\samples^{\Z,(b)}, \samples^{\Dict,(b)}) \right)
\end{equation}
where \bk{$\lossweight_\X=1$ and $\lossweight_\Z \in \R$} may be used to weight the loss terms.
Optimisation computes gradients through both $\samples^{\Z,(i)}$ and $\samples^{\Dict,(i)}$.

The observables $\Dict_{\param_{\Dict}}$ have no natural scale so a \bk{further} regularisation term is used to prevent the observables from collapsing to zero or exploding to infinity.
Recall that the observation variables are considered pairwise as defining complex-valued observations \eqref{eq:background.realrepresentation}.
Regularisation will push the magnitude of each complex pair towards one to encourage \bmk{unit} amplitude,
\begin{equation} \label{eq:regularisation}
  \regularisation_{\Z}(\param; \batch)
  = \frac{1}{\nBatch(\bmk{T}+1) \nSamples \tfrac{\nZ}{2}} \sum_{b=1}^{\nBatch}\sum_{k=0}^\bmk{T} \sum_{j=1}^\nSamples\sum_{\iZ=1,3,\dots,\nZ-1} \left(
      \lvert \gt{z}_{k,\iZ}^{(b,j)} + \im \gt{z}_{k,\iZ+1}^{(b,j)} \rvert - 1
    \right)^2.
\end{equation}
\bmk{The regularisation} will be weighted by a small \bk{factor, see \cref{tab:results.curricula} for the values used in practice.}
The loss \eqref{eq:loss} and the regularisation term \eqref{eq:regularisation} are combined to give the final optimisation target,
\begin{equation}
\label{eq:loss.final}
 \loss(\param; \batch) := \loss_{\X,\Z}(\param; \batch) + \lossweight_{\regularisation}\regularisation_{\Z}(\param; \batch).
\end{equation}

To aid training, a curriculum is used with \bmk{variable} time horizon $\tau$ replacing $T$ in \eqref{eq:loss.state} through \eqref{eq:regularisation}.
In addition to the time horizon $\tau$, each curriculum stage defines loss weights $\lossweight_\Z$, $\lossweight_{\regularisation}$, sample count $\nSamples$, initial and final learning rates $\lr$ and $\lr_{\min}$ for a cosine scheduler and the number of optimiser steps $\nSteps$.

By slowly growing the time horizon $\tau$, even if the spectral radius $\rho$ for dynamics matrix $\K$ is greater than one, then $\Vert \z_\tau \Vert_2 \leq \rho^{\tau} \Vert \z_0 \Vert_2 < \rho^T \Vert \z_0 \Vert_2$, and optimising through the rollout of length $\tau$ is less likely to cause numerical issues than the full $T$ step rollout, particularly when $\tau \ll T$.
In late training, for stable dynamics the optimised linear dynamics $\K$ should have a spectral radius that is at most one (up to some approximation error), \bk{that is} $\rho \leq 1 + \epsilon$ for some small $\epsilon \in \R$.
Then, it is safer to perform the full rollout without numerical issues as $\Vert \z_T\Vert_2 \leq \rho^T \Vert \z_0 \Vert_2 \leq (1+\epsilon)^T \Vert \z_0 \Vert_2 \approx (1 + T\epsilon) \Vert \z_0 \Vert_2$ for $ T \epsilon \ll 1$.

The loss \eqref{eq:loss} sums over all timesteps $k=0,1,\dots,\tau$ for each curriculum stage, \bmk{and the loss }will generally grow as the training passes through stages, approximately by the ratio of time horizons $\tfrac{\tau_{\text{next stage}}}{\tau_{\text{current stage}}}$.

\section{Numerical Experiments}
\label{sec:results}
The probabilistic forecasting approach described in \Cref{sec:method} is now applied to a series of common test problems from the literature.
\bk{Experiments are chosen to assess:
\begin{itemize}
    \item the necessity of probabilistic forecasting (in particular, in \Cref{sec:results.duffing.prediction}),
    \item the ability of a probabilistic model to produce well-calibrated predictive forecasts (via CRPS and spread--skill ratio),
    \item the identification of known eigenfunctions (in particular, \Cref{sec:results.analytic}),
    \item the identification of known eigenvalue frequencies (in particular, in \Cref{sec:results.vortex,sec:results.cavity}),
    \item comparison with a deterministic approach (all) and deep ensemble methods (\Cref{sec:results.analytic,sec:results.duffing}).
\end{itemize}
}

\bk{The test cases are ordered by increasing dataset complexity.}
The first test case is a two dimensional non-linear ODE system with known analytic Koopman eigenfunctions that has been previously studied in, e.g.,~\cite{Lusch2018,Brunton2022,Pan2020}.
This experiment demonstrates our model's ability to identify known eigenfunctions and linear dynamics, and reconstruct the state accurately.
The second test case is the widely studied Duffing oscillator~\cite{Williams2015,Takeishi2017b,Otto2019,Pan2020,Otto2024}.
Again this uses a two dimensional state space, but the observables must now separate two distinct basins of attraction.
\bk{In particular, the need for probabilistic forecasts is highlighted {by an example using noisy data and an} initial condition at the unstable fixed point of the system, in which either stable fixed point could potentially be attained.}
Two examples in fluid dynamics are studied.
Firstly, the {K\'arm\'an} vortex {shedding} phenomenon~\cite{Otto2019,Pan2020,Erichson2019,Colbrook2024} is approximated to indicate how the method extends to high dimensional state spaces (i.e., where the state is a vector of values from a spatial discretisation).
\bk{Then,} a test case investigating cavity flow is studied~\cite{Arbabi2017}, demonstrating the ability to capture a lattice structure of known Koopman-mode frequencies.

Each test problem is studied using multiple dataset configurations: a baseline noise-free dataset (labelled 0), and smaller \bk{datasets} corrupted by additive Gaussian noise $\Ndist(\vec{0},\sigma^2 I_\nX)$,
\begin{equation}
   \gt{\x} + \vec{e}, \vec{e} \sim \Ndist(\vec{0},\sigma^2 I_{\nX}) \text{ for each state }\gt\x \in \R^\nX.
\end{equation}
\bk{
In practice, it is common to have limited training data, for example if the underlying \bmk{data generation is} expensive.
Reducing the dataset size may increase the uncertainty in all three parts of the predictive chain, that is observation, dynamics, and reconstruction, as less of the state space is visited in training.
Additive Gaussian noise is used to represent measurement uncertainty that is often present in real-world datasets.
}
\bk{Models are evaluated against test datasets with the same additive noise to assess the ability to capture state data uncertainty.}

The configurations are summarised in \Cref{tab:results.datasets} and are referenced by their labels throughout this section.
Experiment \bk{D3} is for exposition of multi-modal forecasting so requires only one test trajectory initialised at the origin.
\begin{table}
\footnotesize
\centering
\caption{Dataset configurations for the numerical experiments.
$\nData$ denotes the number of training trajectories, $\sigma$ the additive Gaussian noise standard deviation, $T_{\text{train}}$ the training trajectory length in steps, $\nTest$ and $T_{\text{test}}$ the number and length of the test trajectories, and $\nSamples$ the number of \bk{forecast} samples used in testing.}
\label{tab:results.datasets}
\begin{tabular}{llrcrrrr}
\toprule
Problem & Label & $\nData$ & $\sigma$ & $T_{\text{train}}$ & $\nTest$ & $T_{\text{test}}$ & $\nSamples$\\
\midrule
\multirow{3}{*}{Non-linear ODE}
  & A0 & $5000$  & --        & $50$  & $100$ & $1000$ & $1000$\\
  & A1 & $250$   & $10^{-2}$ & $50$  & $100$ & $1000$ & $1000$\\
  & A2 & $250$   & $10^{-1}$ & $50$  & $100$ & $1000$ & $1000$\\
\midrule
\multirow{4}{*}{Duffing oscillator}
  & D0 & $1500$  & --        & $50$  & $100$ & $200$  & $500$\\
  & D1 & $250$   & $10^{-2}$ & $50$  & $100$ & $200$  & $500$\\
  & D2 & $250$   & $10^{-1}$ & $50$  & $100$ & $200$  & $500$\\
  & D3 & $10000$ & $10^{-2}$ & $800$ & $1$   & $800$  & $500$\\
\midrule
\multirow{3}{*}{Vortex shedding}
  & V0 & $50$    & --        & $50$  & $10$  & $250$  & $20$\\
  & V1 & $10$    & $10^{-2}$ & $50$  & $10$  & $250$  & $20$\\
  & V2 & $10$    & $10^{-1}$ & $50$  & $10$  & $250$  & $20$\\
\midrule
\multirow{3}{*}{Cavity flow}
  & C0 & \bk{$50$}    & --        & $50$  & $10$  & $1000$ & $40$\\
  & C1 & \bk{$10$}    & $10^{-2}$ & $50$  & $10$  & $1000$ & $40$\\
  & C2 & $10$    & $10^{-1}$ & $50$  & $10$  & $1000$ & $40$\\
\bottomrule
\end{tabular}

\end{table}
For each experiment, a deterministic model (denoted \bk{`ML-KMD')} is trained alongside the probabilistic model \bk{(`Prob-KMD')} by removing the noise injection and removing the spread term from the CRPS loss \eqref{eq:background.crps.estimator} to leave only the absolute error.
The deterministic approximation is representative of {machine learning} KMD models, e.g.,~\cite{Lusch2018,Yeung2019,Erichson2019}, noting that we use absolute error compared to mean-square error in the existing literature.
\bk{For the PDE problems, a number of the ML-KMD training instances collapsed to the mean state (i.e., constant state with all dynamics eigenvalues equal to one) and had to be trained under different seeds.} 
\bk{Eigenvalue collapse was not observed in our Prob-KMD runs.}

A direct comparison with the probabilistic methods of~\cite{Pan2020,Otto2024} is not performed \bk{given the challenges in designing and implementing a fair test due to the differences in the model architecture and training:~\cite{Pan2020,Otto2024} fit parametric distributions via likelihood-based inference whereas Prob-KMD fits a sample-based model through a scoring rule.}
\bk{Instead, for the ODE test problems (\Cref{sec:results.analytic,sec:results.duffing}), a comparison with \emph{deep ensembles}~\cite{Lakshminarayanan2017} is included and these will be labelled DE-KMD.
This is an ensemble of models using deterministic observations and dynamics, and reconstruction with a Gaussian head (learned mean and standard deviation). 
These results are generated by training $\nEnsemble=5$ models with different initial seeds, using an identical training strategy.
\bmk{This comparison} favours the deep ensembles which use approximately $\nEnsemble=5$ times the training budget of Prob-KMD\bmk{, and due to this additional computational cost DE-KMD is only applied to the simpler ODE examples.}
}

\bk{The models are compared by computing the mean square error of the ensemble mean, CRPS, and spread--skill ratio.}
\bk{Medians and quartiles are taken over the $\nTest$ test trajectories at each timestep $k=0,1,\dots,T_{\text{test}}$.}
\bk{The ensemble mean is used in the MSE to give a fair comparison between the probabilistic and deterministic model --- otherwise the probabilistic model is penalised for its spread.}

\bk{For a set of forecast samples $\samples$, and a target} \bk{datum} \bk{$\omega$, SSR is defined as~\cite{Fortin2014},}
\begin{equation}\label{eq:background.ssr.single}
    \bk{\ssr(\samples, \omega) :=
    \sqrt{\frac{\nSamples+1}{\nSamples}\frac{\frac{1}{\nSamples-1} \sum_{j=1}^{\nSamples} (x^{(j)} - \bar{x})^2}
    {(\omega - \bar{x})^2}}} = \frac{\text{Spread}}{\text{Skill}},
\end{equation}
\bk{where $\bar{x} = \frac{1}{\nSamples} \sum_{j=1}^{\nSamples} x^{(j)}$ is the sample mean and $\omega$ is the target} \bk{datum}.
\bk{At each forecast step $k$, SSR is aggregated by computing each of the spread and skill as a sum over all state components and test trajectories,}
\begin{equation}\label{eq:background.ssr}
    \bk{\ssr(k) := \sqrt{\frac{\nSamples+1}{\nSamples} \frac{\sum_{b=1}^{\nTest}\sum_{j=1}^{\nX} \frac{1}{\nSamples-1}\sum_{i=1}^{\nSamples} \left(x_{k,j}^{(b,i)} - \bar{x}_{k,j}^{(b)}\right)^2}{\sum_{b=1}^{\nTest}\sum_{j=1}^{\nX} \left(\omega_{k,j}^{(b)} - \bar{x}_{k,j}^{(b)}\right)^2}},}
\end{equation}
\bk{where $x_{k,j}^{(b,i)}$ denotes the $i$-th forecast sample for component $j$ of test trajectory $b$ at step $k$, $\bar{x}_{k,j}^{(b)} := \frac{1}{\nSamples}\sum_{i=1}^{\nSamples} x_{k,j}^{(b,i)}$ is the corresponding sample mean, and $\omega_{k,j}^{(b)}$ is the ground truth target value.}
\bk{A well-calibrated forecast should have $\ssr \approx 1$.}

The ODE problems use multi-layer perceptrons for observation and reconstruction, as outlined in \Cref{sec:appendix.mlp}.
The PDE problems use convolutional neural networks, described in \Cref{sec:appendix.cnn}.
The architectures are outlined in \Cref{tab:results.architectures}.
The training curricula are listed in \Cref{tab:results.curricula}.

\subsection{Non-linear Ordinary Differential Equation with Known Koopman Eigenfunctions}
\label{sec:results.analytic}
{This experiment considers a simple ODE problem with known analytic Koopman eigenfunctions \bk{and eigenvalues,} \bk{thereby facilitating direct comparison with learnt observables $\Dict_{\param_\Dict}$ and eigenvalues of $\K_{\param_\K}$.
The experiment will also demonstrate the predictive uncertainty output by the probabilistic model, and how this evolves through time.}}

\paragraph{Setup}
Consider the state space $\X=[-0.5,0.5]^2$ and the dynamics given by the ordinary differential equation
\begin{equation}\label{eq:analytic_nonlinear.ode}
    \begin{cases}
        \frac{d}{dt} x_1(t) = \alpha x_1(t), \\
        \frac{d}{dt} x_2(t) = \beta (x_2(t) - x_1(t)^2),\\
    \end{cases}
    \text{ for }t>0, \quad \x(0) = \x_0 \in \X,
\end{equation}
for $\alpha=-0.05$ and $\beta=-1$.
Equation \eqref{eq:analytic_nonlinear.ode} is a non-linear ODE with stable fixed point $\x_*=\vec{0}$ and known analytic Koopman eigenfunctions~\cite{Brunton2022},
\begin{equation}\label{eq:analytic_nonlinear.formula}
{\eigfn_1(\x)=}x_2 - \tfrac{\beta}{\beta - 2 \alpha} x_1^2, \quad {\eigfn_2(\x)= }x_1.
\end{equation}
{The vector field $\frac{d \x}{dt}(\x)$ is shown in \Cref{fig:analytic_nonlinear.vectorfield}.
Trajectories on intervals $[0,1]$ are shown in \Cref{fig:analytic_nonlinear.trajectories} for forty different initial conditions $\x_0$ sampled from $\X$.
\begin{figure}\centering
    \hfill
    \subfloat[Vector field $\frac{d \x}{dt}(\x)$]{\label{fig:analytic_nonlinear.vectorfield}%
        \includegraphics{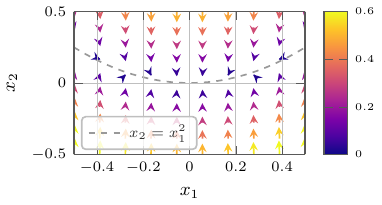}}
    \hfill
    \subfloat[$\x(t)$ on ${[0,1]}$ for forty initial conditions]{\label{fig:analytic_nonlinear.trajectories}%
        \includegraphics{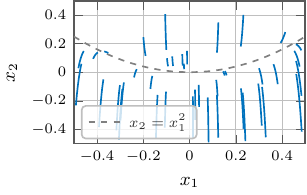}}
    \hfill
    \caption{Dynamics of the non-linear ODE \eqref{eq:analytic_nonlinear.ode}.}
    \label{fig:analytic_nonlinear.dynamics}
\end{figure}}

\paragraph{Data generation}
Discrete time data samples are generated using the \python{} \texttt{solve\_ivp} numerical ODE solver with tolerances $\texttt{rtol}=10^{-8}$ and $\texttt{atol}=10^{-10}$ and the approximation output at intervals of $\Delta t = 0.02$ over the interval $[0,1]$, {resulting in trajectory length} $T_{\text{train}}=50$.
Initial conditions are generated via Latin hypercube sampling over the state space for a given dataset size $\nData$.
For the experiment A0, we choose dataset size $\nData=5000$ to mimic the results in~\cite{Lusch2018}.
We also consider \bk{experiments A1} \bk{and A2} using only $\nData=250$ training trajectories, with the training data corrupted by additive Gaussian noise with standard \bk{deviations $\sigma=10^{-2}$} \bk{and $\sigma=10^{-1}$, respectively}.

\begin{figure}\centering
    \includegraphics{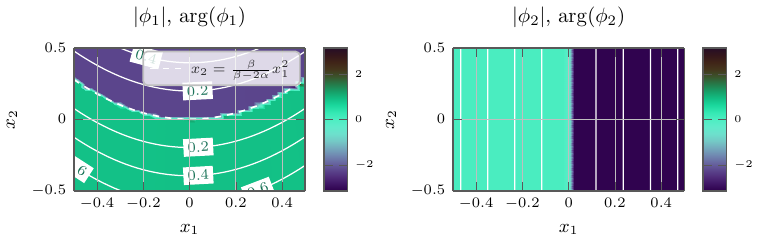}
    \caption{Magnitude (contours) and phase (colour) \bk{of} the expected learnt \bk{observables $\E[\eigfn_1(\x)]$ and $\E[\eigfn_2(\x)]$} for the non-linear ODE experiment \bk{A0 (no noise).}}
    \label{fig:analytic_nonlinear.observables}
\end{figure}

\begin{figure}\centering
    \includegraphics{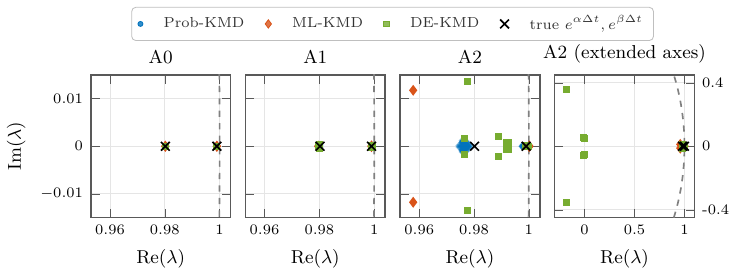}
\caption{\bk{Eigenvalues for the non-linear ODE experiments A0 (no noise), A1 ($\sigma=10^{-2}$) and A2 ($\sigma=10^{-1}$).}}
\label{fig:analytic_nonlinear.dynamics.eigenvalues}
\end{figure}

\paragraph{Observables}
For the noise-free experiment A0, the trained model identifies observables $\Dict_{\param_{\Dict}}(\x)$ that approximate the eigenfunctions of the Koopman operator for the system in \eqref{eq:analytic_nonlinear.ode}.
\bmk{The approximate eigenfunctions are} illustrated by the magnitude and phase of the expected value of the complex-valued functions
\begin{equation}\label{eq:analytic_nonlinear.observable}
\begin{aligned}
\eigfn_1(\x) := \Dict_{\param_{\Dict},1}(\x) + \im\Dict_{\param_{\Dict},2}(\x),\\
\eigfn_2(\x) := \Dict_{\param_{\Dict},3}(\x) + \im\Dict_{\param_{\Dict},4}(\x),
\end{aligned}
\end{equation}
which are shown in \Cref{fig:analytic_nonlinear.observables}.
These approximate (up to proportionality constant) the true \bk{eigenfunctions given in \eqref{eq:analytic_nonlinear.formula}.} %
\bk{The discrete-time eigenvalues are shown in \Cref{fig:analytic_nonlinear.dynamics.eigenvalues}.}
\bk{In experiments A0 and A1, the model eigenvalues approximate the corresponding true eigenvalues $0.999\approx \exp(\alpha \Delta t)$ and $0.980 \approx \exp(\beta \Delta t)$.}
\bk{It is only in the experiment A2, with the largest additive noise, that the models differ.}
\bk{Prob-KMD eigenvalues remain accurate.}
\bk{Three of the five ensemble members of DE-KMD do not recover the `fast' (furthest from modulus one, $0.980$) behaviour, appearing to rely on the Gaussian head to instead capture this as state uncertainty.}
\bk{ML-KMD also does not capture the fast eigenvalue correctly, estimating approximately $0.958\pm0.012 \im$ instead of the true $0.980$.}

\paragraph{Predictive performance}\label{sec:analytic_nonlinear.prediction}
Model performance is analysed using $\nTest=100$ initial conditions sampled from the state space $\X$ via Latin hypercube sampling.
Trajectory data is generated on an extended interval $[0,20]$ to give $1000$ steps.

\begin{figure}\centering
    \includegraphics{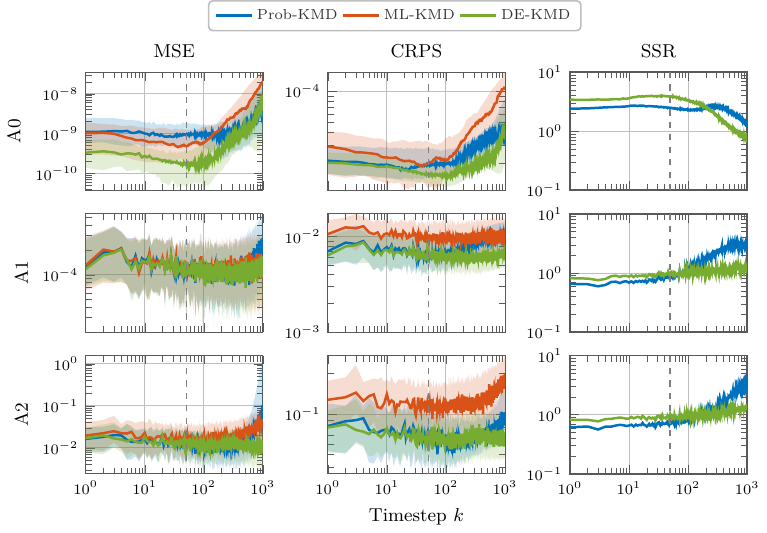}
    \caption{Median and quartile \bk{MSE} and CRPS, and \bk{SSR, against step $k$} for the non-linear ODE experiments \bk{A0 (no noise), A1 ($\sigma=10^{-2}$)}\bk{ and A2 ($\sigma=10^{-1}$)}\bk{. Dashed line: end of training horizon.}}
    \label{fig:analytic_nonlinear.prediction}
\end{figure}

The A0 experiment trains accurate models with mean-square error that remains small even when predicting to $T_{\text{test}}=1000$ steps, as shown in \bk{\Cref{fig:analytic_nonlinear.prediction}.}
The A1 experiment has mean-square error of approximately $10^{-4}$, which is the limit imposed by the additive noise.
\bk{Similarly the A2 experiment sits at approximately $10^{-2}$.}
\bk{At $k=1000$, the Prob-KMD model has the largest error in both experiments A1 and A2, but the models have mostly comparable performance.}

\begin{figure}\centering
    \subfloat[\bk{Prob-KMD}]{\includegraphics{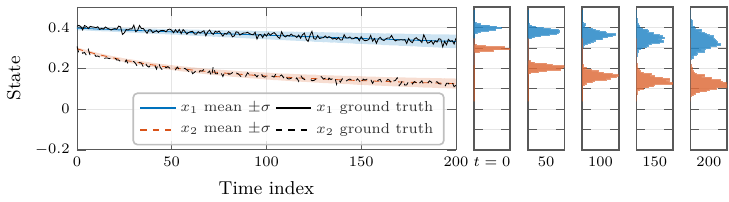}\label{fig:analytic_nonlinear.pdfs.noisy.prob}}\\
    \subfloat[\bk{DE-KMD}]{\includegraphics{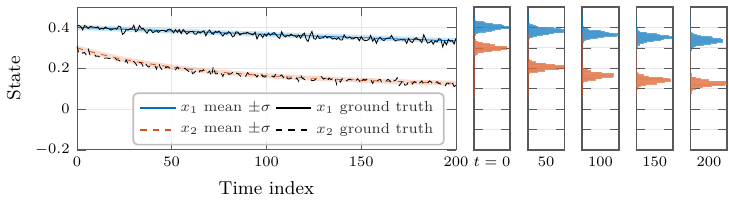}\label{fig:analytic_nonlinear.pdfs.noisy.ens}}
    \caption{Example trajectory \bk{for the non-linear ODE experiment A1 ($\sigma=10^{-2}$).} \bk{Expected state $\E[x_k]$ $\pm$ one standard deviation}, ground truth (black) and snapshot marginal \bk{PDFs.}}
    \label{fig:analytic_nonlinear.pdfs.noisy}
\end{figure}

The novelty \bk{of Prob-KMD is its CRPS-optimised probabilistic forecast.}
\bk{The probabilistic models are seen to outperform the deterministic models in CRPS, with the DE-KMD model slightly outperforming the Prob-KMD model, particularly on long rollouts.}
Snapshots of the marginal PDFs in \Cref{fig:analytic_nonlinear.pdfs.noisy} show the A1 \bk{probabilistic models capture} the additive Gaussian noise that is present in the training data.
\bk{The over-dispersion of Prob-KMD becomes visible for large $k$.}
The spread--skill ratio, computed as \eqref{eq:background.ssr}, and shown in \bk{\Cref{fig:analytic_nonlinear.prediction}}, also highlights \bk{the over-dispersion.}
\bk{The A1 Prob-KMD model is slightly under-dispersed over the training horizon.}
\bk{When rolling out to $T_{\text{test}}=1000$, twenty times longer than the training data, the model gradually becomes over-dispersed ($\ssr_{1000} \approx 3 > 1$).}
The \bk{A1 DE-KMD model has a spread--skill ratio close to one for all times.}
\bk{Similar performance is seen for both models in the A2 experiment.}
\bk{In the} A0 experiment \bk{both models are over-dispersed, with DE-KMD more so.}
\bk{This is hypothesised to be} a consequence of the \bk{high model skill, and the model not suppressing the residual spread enough.}

\subsection{Duffing Oscillator}
\label{sec:results.duffing}
\begin{figure}\centering
    \hfill
    \subfloat[Vector field $\frac{d \x}{dt}(\x)$]{\label{fig:duffing.vectorfield}%
        \includegraphics{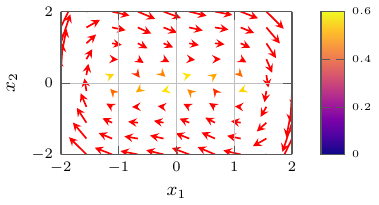}}
    \hfill
    \subfloat[$\x(t)$ on ${[0,1]}$ for forty initial conditions]{\label{fig:duffing.trajectories}%
        \includegraphics{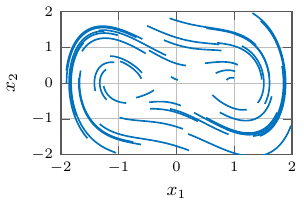}}
    \hfill
    \caption{Dynamics of the Duffing oscillator \eqref{eq:duffing.ode}.}
    \label{fig:duffing.dynamics}
\end{figure}

{This problem considers the Duffing oscillator, and in particular \bk{includes an experiment highlighting the challenges from} diverging dynamics.}
\paragraph{Setup}
Consider the state space $\X=[-2,2]^2$ and the dynamics given by the ordinary differential equation
\begin{equation}\label{eq:duffing.ode}
    \begin{cases}
        \frac{d}{dt} x_1(t) = x_2(t) \\
        \frac{d}{dt} x_2(t) = -\delta x_2(t) - \beta x_1(t) - \alpha x_1(t)^3 \\
    \end{cases}
    , \text{ for }t>0, \quad \x(0) = \x_0 \in \X,
\end{equation}
for $\alpha=1$, $\beta=-1$, and $\delta=0.5$.
Equation \eqref{eq:duffing.ode} is a non-linear ODE representing the Duffing oscillator with stable fixed points at $\x_{*,\pm}=[\pm \sqrt{-\beta/\alpha},0]=[\pm 1,0]$ and an unstable fixed point at the origin $\x_{*,0}=\vec{0}$~\cite{Williams2015,Otto2019}.
The vector field $\frac{d \x}{dt}(\x)$ is illustrated in \Cref{fig:duffing.vectorfield}.
Trajectories on intervals $[0,1]$ are shown in \Cref{fig:duffing.trajectories} for forty different initial conditions $\x_0$ sampled from $\X$.

{At the stable fixed points $\x_{*,\pm}$, the continuous-time dynamics \eqref{eq:duffing.ode} have Jacobian
\begin{equation}\label{eq:duffing.jacobian}
    \begin{bmatrix}
        0 & 1\\
        2\beta & -\delta
    \end{bmatrix}
\end{equation}
with eigenvalues $-0.5{\delta} \pm 0.5\sqrt{\delta^2 + 8 \beta} \approx -0.25 \pm \im 1.392${, noting that $\delta^2 + 8 \beta < 0$}.
Converting to the discrete-time equivalent for a timestep of $\Delta t = 0.02$, the linearised discrete-time dynamics have eigenvalues
\begin{equation}\label{eq:duffing.eigenvalues}
\begin{aligned}
    \lambda_{*,\pm} &= \exp(-0.5{\delta} \Delta t)\exp(\pm \im 0.5\sqrt{-\delta^2 - 8 \beta} \Delta t) \\
    &\approx 0.9950 \exp(\pm \im 0.0278) \\
    &\approx 0.9946 \pm \im  0.0277.
\end{aligned}
    \end{equation}
This linearisation will help us interpret the learnt observables in our probabilistic model.}

\bk{It was observed in preliminary experiments that the model may learn eigenvalues capturing the inverse of the rate of attraction, ${\lvert\lambda_{*,\pm}\rvert^{-1}>1}$, which are equally valid in our setup.
Linear rollout leads to exponential growth instead of convergence to zero, resulting in unbounded predicted observations that are not represented in the training data.
The dynamics are therefore constrained to be stable as described in \Cref{rem:method.spectralradius}, so the temporal rollout is contractive.}

\paragraph{Data generation}
The data generation process is identical to \bk{that} in \Cref{sec:results.analytic}. 
The experiment D0 will use $\nData=1500$ trajectory samples, whilst the noisy, low data \bk{experiments D1} \bk{and D2} will use additive Gaussian noise with standard \bk{deviations $\sigma=10^{-2}$} \bk{and $\sigma=10^{-1}$, respectively,} and $\nData=250$ trajectories.

\begin{figure}\centering
   \includegraphics{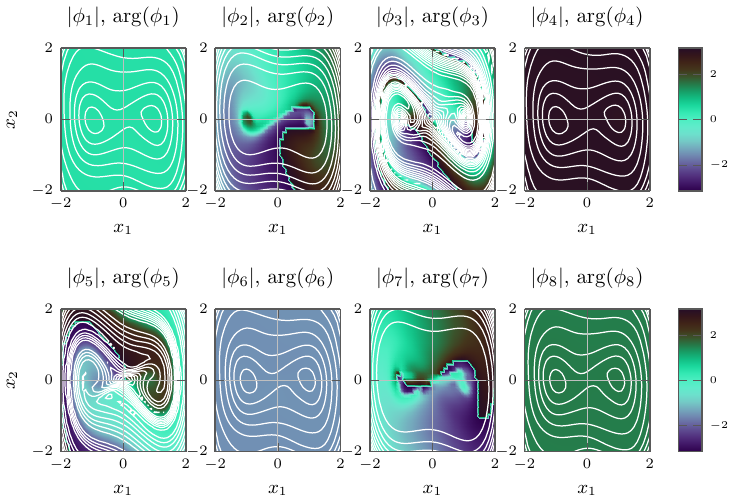}
    \caption{Magnitude (contours) and phase (colour) \bk{of} the expected learnt \bk{observables $\E[\eigfn_j(\x)]$} for the Duffing oscillator experiment \bk{D0 (no noise).}}
    \label{fig:duffing_oscillator.observables}
\end{figure}

\begin{figure}\centering
    \includegraphics{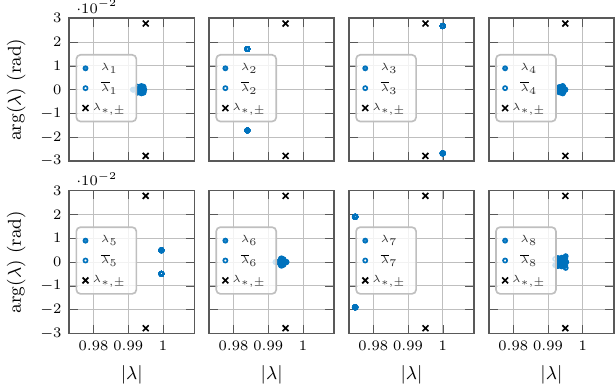}
    \caption{Magnitude and phase \bk{of} samples of the eigenvalues $\{\lambda_i\}_{i=1}^{\tfrac{\nZ}{2}}$ for the Duffing oscillator experiment \bk{D0 (no noise).}}
    \label{fig:duffing.eigenvalues}
\end{figure}

\paragraph{Observation}
The observation functions $\Dict_{\param_{\Dict}}(\x)$ are analysed again as complex-valued functions
\begin{equation}\label{eq:duffing.observable}
    \eigfn_j(\x)=\Dict_{\param_{\Dict},2j-1}(\x) + \im\Dict_{\param_{\Dict},2j}(\x)
\end{equation}
for observable pairs indexed $j=1,2,\dots,\tfrac{\nZ}{2}$.
The magnitude and the phase of the expected value of the observations learnt in experiment D0 are shown in \Cref{fig:duffing_oscillator.observables}.
In particular, the distinct basins of attraction are clearly captured in the complex observables {$\eigfn_3,\eigfn_5$}.
{Samples of the corresponding eigenvalues $\{\lambda_i\}_{i=1}^{{\nZ}/{2}}$ are shown in \Cref{fig:duffing.eigenvalues}.
The eigenvalues $\lambda_{*,\pm}$ corresponding to the \bk{linearised dynamics at the} stable fixed points $\x_{*,\pm}$ are also shown.
\bk{The eigenvalue $\lambda_3$ captures the rotational dynamics around the fixed points.
The eigenvalues $\lambda_1,\lambda_4,\lambda_6,\lambda_8$ capture the exponential decay towards the fixed points.
The eigenvalues $\lambda_2,\lambda_7$ are hard to clearly interpret, but represent the long range rotation and decay in the state space.
Lastly, the eigenvalue $\lambda_5$ is less interpretable.}
\bk{It is interesting that the rotation and decay effects have been separated by the model.}
}

\begin{figure}\centering
    \includegraphics{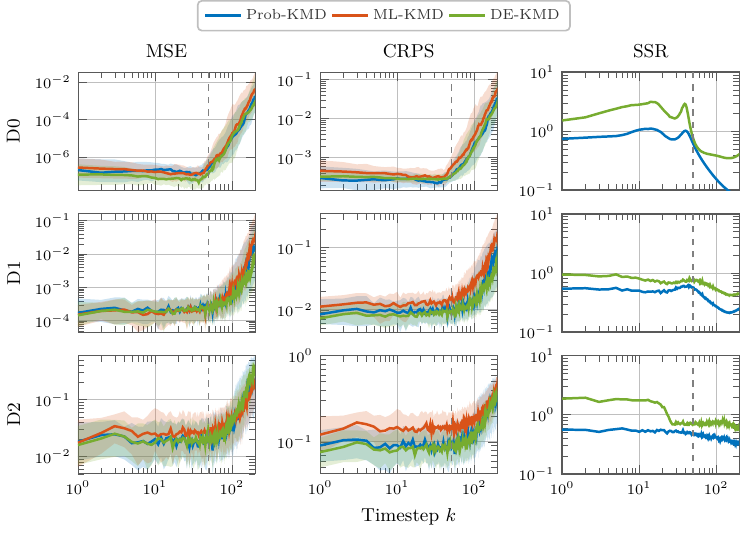}
    \caption{Median and quartile \bk{MSE} and CRPS, and \bk{SSR, against step $k$} for the Duffing oscillator experiments D0 \bk{(no noise), D1 ($\sigma=10^{-2}$)}\bk{ and D2 ($\sigma=10^{-1}$)}\bk{. Dashed line: end of training horizon.}}
    \label{fig:duffing.loss}
\end{figure}

\paragraph{\bk{Predictive} performance}
\label{sec:results.duffing.prediction}
The mean-square prediction error, shown in \bk{\Cref{fig:duffing.loss},} is small over the training \bk{interval, generally comparable to the noise $\sigma^2$ in experiments D1 and D2.}
Beyond the training horizon, the approximation errors increase steadily with step $k$.

The probabilistic \bk{models outperform ML-KMD in CRPS, and Prob-KMD is competitive with the deep ensemble model.}
\bk{Again, for long temporal rollouts there is a steady increase in CRPS in all models.}

The spread--skill ratio is \bk{also} shown in \bk{\Cref{fig:duffing.loss}.}
\bk{In D0, the Prob-KMD model is approximately calibrated over the training interval, and the deep ensemble is initially over-dispersed.}
\bk{Both become increasingly under-dispersed beyond the training interval, in line with the increasing approximation error.}
\bk{Adding noise in the D1 and D2 experiments results in an under-dispersed Prob-KMD model, with some uncertainty captured.}
\bk{DE-KMD is better calibrated in both D1 and D2, although also under-dispersed for long rollouts.}

\paragraph{The need for probabilistic modelling}
\bk{In the} final Duffing oscillator experiment, \bk{D3, we} will demonstrate the importance of the probabilistic model in capturing diverging dynamics in the Duffing oscillator.
\bk{Near the boundary of the basins of attraction, the future state is extremely sensitive to the initial condition.
    It is hard for a data-driven model to perfectly capture the boundary of the basins of attraction, particularly under measurement uncertainty that may be present in an experimental dataset.
    We wish to capture this uncertainty and produce an appropriate multi-modal forecast.}

The additive Gaussian noise is again $\sigma=10^{-2}$, but the \bk{models are} trained with $\nData=10000$ trajectories of length $T_{\text{train}}=800$.
The purpose is to create an accurate model of the dynamics at the unstable fixed point \bk{$\x_{*,0}=\vec{0}$.}
\bk{An extended training curriculum is used to account for the longer trajectories as described in \Cref{tab:results.curricula}.}

\begin{figure}
    \centering
    \subfloat[\bk{Prob-KMD}]{\includegraphics{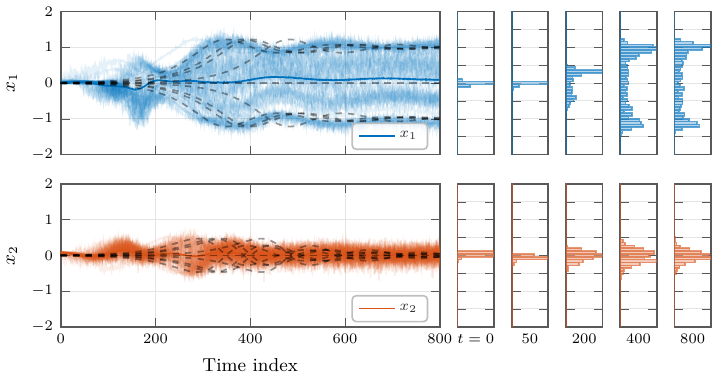}\label{fig:duffing.diverging.prob}}\\
    \subfloat[\bk{DE-KMD}]{\includegraphics{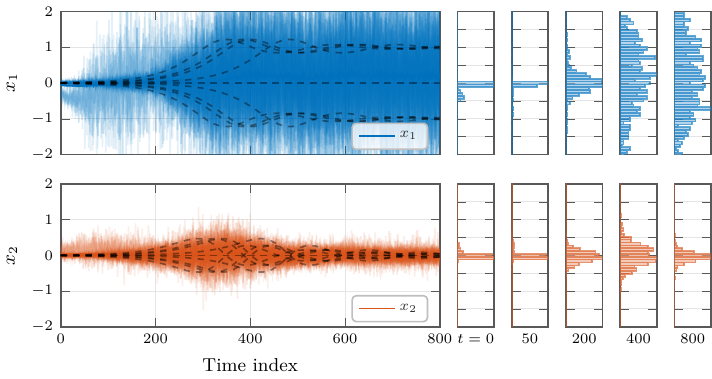}\label{fig:duffing.diverging.ens}}
    \caption{\bk{Trajectories from $\gt\x_0=[0,0]$ for the Duffing oscillator experiment D3 ($\sigma=10^{-2}$):} expected state (solid), forty samples and snapshot marginal \bk{PDFs.}
    Black dashed \bk{lines:} ground truth \bk{without noise from} initial conditions $\{-0.01,0,0.01\}^2$.}
    \label{fig:duffing.diverging}
\end{figure}

The initial condition is chosen to be $\gt\x_0:=[0,0]$ which is the unstable fixed point of the dynamics \eqref{eq:duffing.ode}.
The trained \bk{Prob-KMD} model captures the uncertainty in the input data, shown in \bk{\Cref{fig:duffing.diverging.prob},} where samples of the reconstructed state after observation $\x_0 \sim g_{\param_{g}}(\Dict_{\param_\Dict}(\gt\x_0))$ take positive and negative values in \bk{both $x_1$ and $x_2$.}
\bk{A} probabilistic model should then recognise that the state could be attracted to either of the stable fixed points.
Approximate trajectories and snapshot marginal state PDFs are shown in \bk{\Cref{fig:duffing.diverging} for Prob-KMD and DE-KMD models.}
\bk{Prob-KMD} samples capture trajectories that are attracted to both fixed points $\x_{*,\pm}=[\pm 1, 0]$.
\bmk{Both fixed points} \bk{cannot} be \bmk{simultaneously} identified by a deterministic model, which shows the necessity of the probabilistic approach.
\bk{The deep ensemble is a mixture of $\nEnsemble=5$ Gaussians, so offers a probabilistic approximation, but is less expressive.}
\bk{Samples of the DE-KMD model can be seen in \Cref{fig:duffing.diverging.ens} to spread to cover both fixed points, but it has not captured the distinct basins of attraction.}
\bk{The forecast is essentially a Gaussian peak with large variance.}
\bk{Prob-KMD is superior at resolving the multi-modal forecast in this experiment.}

\subsection{Vortex Shedding}
\label{sec:results.vortex}
To demonstrate an application to high dimensional state data, the von K\'arm\'an vortex street phenomenon is studied.
Vortex shedding is a popular test problem in the literature~\cite{Rowley2009,Bagheri2013,Colbrook2024} with practical relevance in engineering design.

\begin{figure}\centering
    \subfloat[Expected value]{\label{fig:vortex.reconstruction.mean}%
        \includegraphics[]{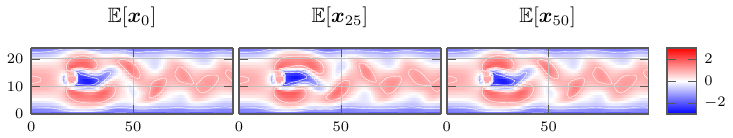}}\\
    \subfloat[Absolute error]{\label{fig:vortex.reconstruction.abserror}%
        \includegraphics[]{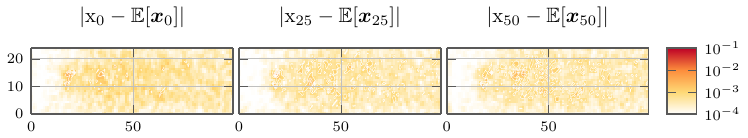}}\\
    \subfloat[Spread--skill ratio]{\label{fig:vortex.reconstruction.ssr}%
        \includegraphics[]{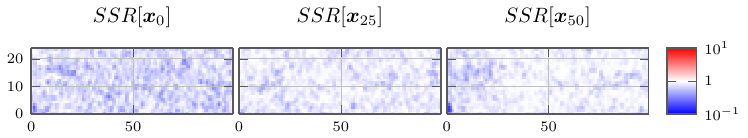}}
    \caption{\bk{Prob-KMD expected state} $\E[\x]$, absolute error $\lvert \gt{\x} - \E[\x] \rvert$ and \bk{pointwise SSR} $\ssr[\x]$ at steps $k=0,25,50$ for \bk{a test trajectory of} the vortex shedding experiment \bk{V0 (no noise).}}
    \label{fig:vortex.reconstruction}
\end{figure}

\paragraph{Setup}
The problem setup is taken from~\cite{Colbrook2024}.
Briefly, the velocity field is computed for a flow in a two-dimensional channel of height $5D$ with a cylinder of diameter $D$ vertically centred, as seen in the \bk{example forecast} states shown in \Cref{fig:vortex.reconstruction}.
A parabolic inflow is applied on the left boundary and the {flow} separates {downstream of} the cylinder.
The Reynolds number $\Reynolds=100$ results in laminar vortex shedding in the cylinder wake~\cite{Colbrook2024}.

Vortex shedding is a popular test case in the DMD literature as the state cycles with {a} shedding frequency $f_0$.
{The flow} can be well approximated using linear projection and linear reconstruction via a sequence of increasingly oscillatory spatial modes with temporal angular frequencies $k f_0$ for $k \in \N$~\cite{Bagheri2013,Colbrook2024}.

\paragraph{Data generation}
\bk{We} use the horizontal velocity $u_x$ \bk{from~\cite{Colbrook2024}} on a spatial grid of $200$ vertical cells by $800$ horizontal cells from step $k=0$ to \bk{$1000$,} with approximately $48$ steps per shedding cycle.
To reduce the state space size, the spatial grid is downsampled by a factor of eight {in both dimensions}, to give a state space $\X=\R^{25 \times 100}$.
\bmk{Downsampling} allows us to reduce the model capacity and training time, without drastically changing the underlying dynamics that we wish to capture.
Prior to any processing, the entire dataset is rescaled to be zero mean and unit variance.

The training dataset V0 consists of $\nData=50$ \bk{trajectories,} each of length $T_{\text{train}}=50$ timesteps that are sampled from steps $k=0$ to $500$ in the supplied dataset.
The data from steps $k=501$ to $700$ are used during validation, and the steps from $k=701$ to $1000$ are used as test data to study the performance.
The \bk{datasets V1} \bk{and V2} \bk{use} $\nData=10$ training trajectories with additive Gaussian noise with standard deviations \bk{$\sigma=10^{-2}$ and $\sigma=10^{-1}$, respectively.}

Within the training and test datasets individual trajectories may overlap, which is a consequence of the limited dataset and extracting all trajectories from a single long trajectory.

\begin{figure}\centering
\includegraphics{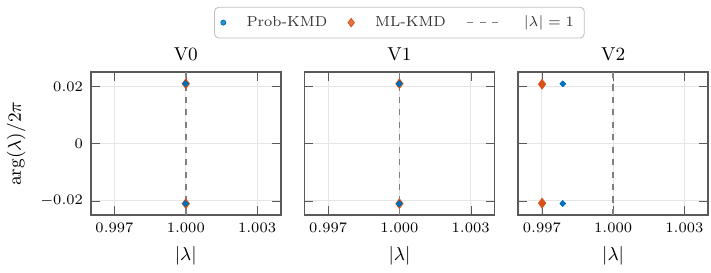}
\caption{\bk{Eigenvalue magnitude} and phase for the vortex shedding experiments V0 \bk{(no noise), V1 ($\sigma=10^{-2}$)}\bk{ and V2 ($\sigma=10^{-1}$)}.}
    \label{fig:vortex.eigenvalues}
\end{figure}

\paragraph{Observation}
To characterise a state $\x\in\X$, the observations simply need to capture the vortex shedding phase.
\bk{Non-linear reconstruction implies a single frequency $f$ can represent the periodic behaviour for all other Koopman mode frequencies $kf$ for $k\in\N$.}
\bk{Linear reconstruction would require the full spectrum lattice.}
To this end, the model is trained with a single pair of observables, $\nZ=2$, which is sufficient to capture the vortex shedding phase.
The corresponding eigenvalues are shown in \bk{\Cref{fig:vortex.eigenvalues}.}
\bk{
The oscillatory behaviour, with period of approximately $48$ timesteps, is well captured by the eigenvalue phase in all three experiments.
This is perhaps not surprising: additive noise will not change the underlying dynamics.
In experiment V2 the eigenvalue moduli are $0.997$ (ML-KMD) and $0.998$ (Prob-KMD), and this dissipation will be seen to introduce forecasting errors for long rollouts.
}

\paragraph{Predictive performance}
Snapshots of the expected value of the predicted trajectories $\E[\x]$, the corresponding absolute error $\lvert \gt{\x} - \E[\x] \rvert$ and the spread--skill ratio $\ssr[\x]$, computed \bk{pointwise at each grid node,} are shown in \Cref{fig:vortex.reconstruction} for a test \bk{trajectory and the Prob-KMD model for experiment V0.}
The absolute error values shown in \Cref{fig:vortex.reconstruction.abserror} demonstrate the high quality approximation.
The spread--skill ratio in \Cref{fig:vortex.reconstruction.ssr} demonstrates the model captures what little uncertainty is present.

\begin{figure}\centering
    \includegraphics{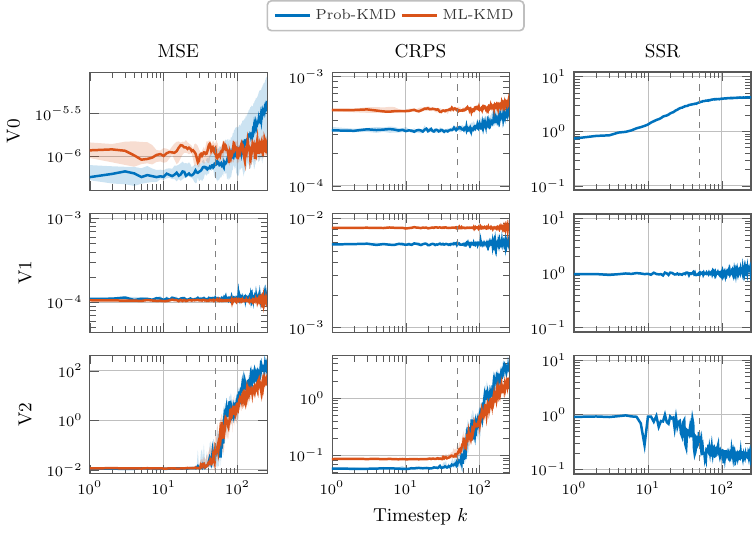}
    \caption{Median and quartile \bk{MSE} and CRPS, and \bk{SSR, against step $k$} for the vortex shedding experiments V0 \bk{(no noise), V1 ($\sigma=10^{-2}$)}\bk{ and V2 ($\sigma=10^{-1}$)}\bk{. Dashed line: end of training horizon.}}
    \label{fig:vortex.loss}
\end{figure}

The test data mean-square error and CRPS for V0 \bk{and V1,} shown in \bk{\Cref{fig:vortex.loss},} are {fairly} stable even when rolled out to a maximum timestep $T_{\text{test}}=250$, five times longer than the training data.
\bk{Prob-KMD outperforms the ML-KMD model in CRPS in experiments V0 and V1.}
\bk{In experiment V2, with $\sigma=10^{-1}$, initially Prob-KMD outperforms ML-KMD up to $k\approx65$, but both grow in error over long rollouts, with ML-KMD having lower CRPS.}
\bk{In V0, the model is briefly well calibrated, before the SSR grows steadily to $\ssr\approx4$.}
\bk{\bmk{High SSR} is hypothesised to be due to the high skill in the noise-free experiment, and residual spread from the noise injection.}
\bk{Additive noise} in experiment V1 \bk{results in a near perfect $\ssr\approx1$ for all timesteps.}

\bk{In V0 and V1 the learned eigenvalues have modulus $1.0000$, so propagated observations remain on the unit circle even with high exponents.}
\bk{In V2, the SSR is initially close to one, before becoming under-dispersed after approximately thirty steps (decreasing to approximately $0.2$ by 100 steps).}
\bk{This diverging forecast is likely a consequence of the exponential decay} due to \bk{the learned eigenvalues having modulus less than one.}
\bk{\bmk{Consequently}, in long rollouts, the observation moves away from the constant modulus circle that observations should map to.}
\bk{Based upon the known physics, constraining the eigenvalue magnitude to exactly one could alleviate this.}

\subsection{Cavity Flow}
\label{sec:results.cavity}
Another popular fluid dynamics test case is lid-driven cavity flow, investigated in the Koopman context in, \bk{e.g.,}~\cite{Arbabi2017,Colbrook2025rigged}.
This test problem gives rise to a lattice structure of Koopman-mode frequencies, and we demonstrate that our model is able to \bk{identify this structure.
The spectral redundancy due to non-linear state reconstruction is also discussed}.

\paragraph{Setup}
\bk{We} use the problem setup in \bk{\cite{Arbabi2017}.}
The problem considers a square spatial domain $\spatialdomain=[-1,1]^2$, with a non-zero velocity boundary condition on the lid $\partial \spatialdomain_{lid} = [-1,1] \times \{1\}$ with $u_{\xspace_1}=\left(1-\xspace_1^2\right)^2$ and $u_{\xspace_2}\equiv0$.
The {quantity} of interest is the vorticity $\omega:\spatialdomain \to \R$,
\begin{equation}\label{eq:vorticity}
    \omega(\vxspace) := - \nabla \cdot \nabla \streamfunction(\vxspace),
\end{equation}
where $\streamfunction:\spatialdomain \to \R$ is the stream function, which is related to the velocity field $\vec{u}:\spatialdomain \to \R^2$ via
\begin{equation}
    \vec{u}(\vxspace) = \begin{bmatrix}\partial_{\xspace_2} \streamfunction(\vxspace), -\partial_{\xspace_1} \streamfunction(\vxspace)\end{bmatrix}^{\top}.
\end{equation}
{Specifically, the} state vector $\x$ is the vorticity on the grid of Chebyshev nodes.
The behaviour of the flow is controlled by the Reynolds number $\Reynolds$.
The dataset selected uses $\Reynolds=16 \times 10^{3}$ which is known to produce a quasi-periodic flow~\cite{Arbabi2017}.

\paragraph{Data generation}
\bk{We use the dataset supplied with~\cite{Arbabi2017}.}
The approximation is \bk{computed} using the stream-function formulation of the incompressible Navier--Stokes equations.
The solver uses a spectral collocation method returning the stream function \bk{$\streamfunction$} on a grid $\Xi = \{\xi_1,\xi_2,\dots,\xi_{65}\}^2 \subset \spatialdomain$ of $65 \times 65$ Chebyshev nodes $\{\xi_1,\xi_2,\dots,\xi_{65}\}$ on the interval $[-1,1]$~\cite{Trefethen2000}.
Further methodological details are given in~\cite{Arbabi2017,Shen1991}.
The vorticity $\omega$ is then computed for the collocation nodes via \eqref{eq:vorticity}.
The nodal values of the vorticity form the state $\x\in\R^{65\times65}$,
\begin{equation}
    \x_{i,j}:= \omega(\xi_i, \xi_j) \text{ for } i,j=1,2,\dots,65. 
\end{equation}
State snapshot data is recorded with timestep \bk{$0.1$.}
Again, the entire dataset is normalised to have mean zero and variance one.

The experiment C0 uses a training dataset with $\nData=50$ trajectories with $T_{\text{train}}=50$, the \bk{experiments C1 and C2 use training datasets} with $\nData=10$, $T_{\text{train}}=50$ and additive Gaussian noise with standard \bk{deviations $\sigma=10^{-2}$ and $\sigma=10^{-1}$ respectively.}

\paragraph{Observation}
The {observables $\Dict$} are best analysed using the {corresponding eigenvalues}.
The results in~\cite{Arbabi2017} identify a discrete spectrum of frequencies specified by $f = k_1 f_1 + k_2 f_2$ for $[f_1,f_2] = [0.155375, 0.096911]$ and $k_1,k_2\in\Z$.
The lattice structure of mode frequencies is discussed in further detail in, e.g.,~\cite{Arbabi2017,Brunton2022}.
The expected value of the eigenvalue for each observable is listed in \Cref{table:cavity.eigenvalues} alongside the corresponding lattice frequencies \bk{(one per complex pair)}.
\begin{table}
\footnotesize
\centering
\caption{\centering Learnt eigenvalues \bk{$\lambda$ for the cavity flow experiment C0 (no noise), with} frequency $\arg(\lambda)/(2\pi\cdot0.1)$, lattice coefficients $(k_1,k_2)$ and \bk{lattice} frequency \bk{$k_1 f_1 + k_2 f_2$.}}
\label{table:cavity.eigenvalues}
\begin{tabular}{ccccc}
\toprule
$\Re[\E[\lambda]]$ & $\Im[\E[\lambda]]$ & $\arg(\E[\lambda])/(0.2\pi)$ & $(k_1,k_2)$ & $k_1 f_1 + k_2 f_2$\\
\midrule
$1.0000$ & $0.0000$ & $0.0000$ & $(0,0)$ & $0.0000$\\
$0.9983$ & $0.0609$ & $0.0970$ & $(0,1)$ & $0.0969$\\
$0.9982$ & $0.0609$ & $0.0970$ & $(0,1)$ & $0.0969$\\
$0.9953$ & $0.0973$ & $0.1550$ & $(1,0)$ & $0.1554$\\
$0.9954$ & $0.0973$ & $0.1551$ & $(1,0)$ & $0.1554$\\
$0.9951$ & $0.0973$ & $0.1551$ & $(1,0)$ & $0.1554$\\
$0.9949$ & $0.0973$ & $0.1552$ & $(1,0)$ & $0.1554$\\
$0.9876$ & $0.1579$ & $0.2523$ & $(1,1)$ & $0.2523$\\
\bottomrule
\end{tabular}
\end{table}
A zero frequency component $\lambda \approx 1$ is identified, corresponding to the time-invariant mean flow.
Three distinct non-zero lattice frequencies are identified in the model trained in experiment C0, namely $f_2$, $f_1$ and $f_1 + f_2$, and notably these frequencies are seen to represent {the} modes with the highest kinetic energy~\cite{Arbabi2017}.
There are repetitions, e.g., $f=\pm 0.1554$ can be identified with four observable pairs, and all eigenvalues have magnitude very close to $1$ which is in agreement with the results and the nature of the quasi-periodic flow.
The model is able to identify and learn appropriate Koopman modes for the cavity flow problem and it will be seen below that these modes are enough for a good quality forecast.
\bk{Again, redundancy in the spectrum is due to the non-linear reconstruction, but some higher order modes are identified as they are still consistent with linear observable dynamics.}
{
The identified frequencies are clearly seen in \bk{\Cref{fig:cavity.eigenvalues} for the C0 Prob-KMD model.}
\begin{figure}\centering
\includegraphics{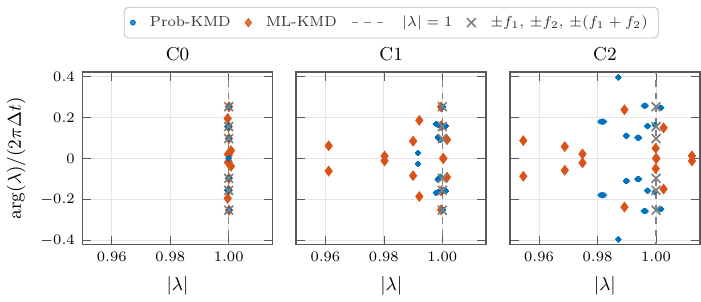}
\caption{\bk{Eigenvalue magnitude and frequency} for the cavity flow experiments C0 \bk{(no noise), C1 ($\sigma=10^{-2}$)}\bk{ and C2 ($\sigma=10^{-1}$)}.}
\label{fig:cavity.eigenvalues}
\end{figure}}
\bk{The C0 ML-KMD model captures these frequencies, some additional redundant lattice frequencies, and spurious additional frequencies.}
\bk{\bmk{Spurious frequencies} become more pronounced in experiments C1 and C2.}
\bk{Whilst the phase is still generally well captured in C1, and to some extent C2, the eigenvalue magnitudes become increasingly far from one as the training dataset noise increases, particularly for the ML-KMD models with many introducing damping in time.}

\begin{figure}\centering
    \subfloat[Expected value]{\label{fig:cavity.reconstruction.mean}%
        \includegraphics[]{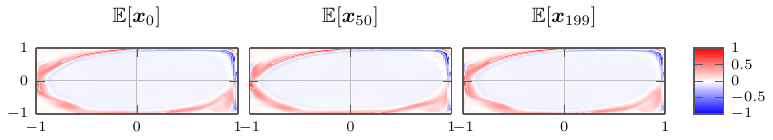}}\\
    \subfloat[Absolute error]{\label{fig:cavity.reconstruction.abserror}%
        \includegraphics[]{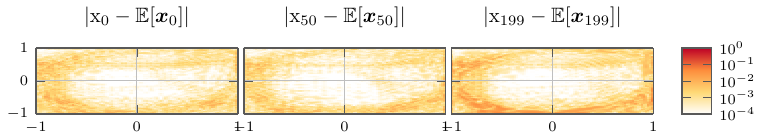}}\\
    \subfloat[Spread--skill ratio]{\label{fig:cavity.reconstruction.ssr}%
        \includegraphics[]{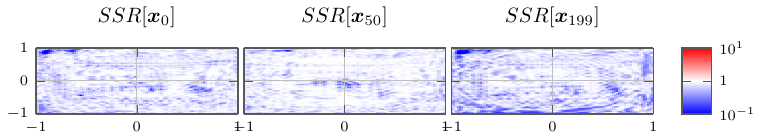}}
    \caption{\bk{Prob-KMD expected state} $\E[\x]$, absolute error $\lvert \gt{\x} - \E[\x] \rvert$ and \bk{pointwise SSR} $\ssr[\x]$ at steps \bk{$k=0,50,199$ for a test trajectory of} the cavity flow experiment \bk{C0 (no noise).}}
    \label{fig:cavity.reconstruction}
\end{figure}

\paragraph{Predictive performance}
The absolute approximation error in \Cref{fig:cavity.reconstruction.abserror} indicates there is a generally small approximation error, concentrated in the main recirculating flow, that grows with timestep.
The \bk{{pointwise}} spread--skill ratio in \Cref{fig:cavity.reconstruction.ssr} indicates that the probabilistic model captures the uncertainty well with $\ssr$ generally close to one.

The mean-square error, CRPS, and SSR are computed for a test set of $\nTest=10$ \bk{trajectories} with length \bk{$T_{\text{test}}=1000$.}
\begin{figure}\centering
    \includegraphics{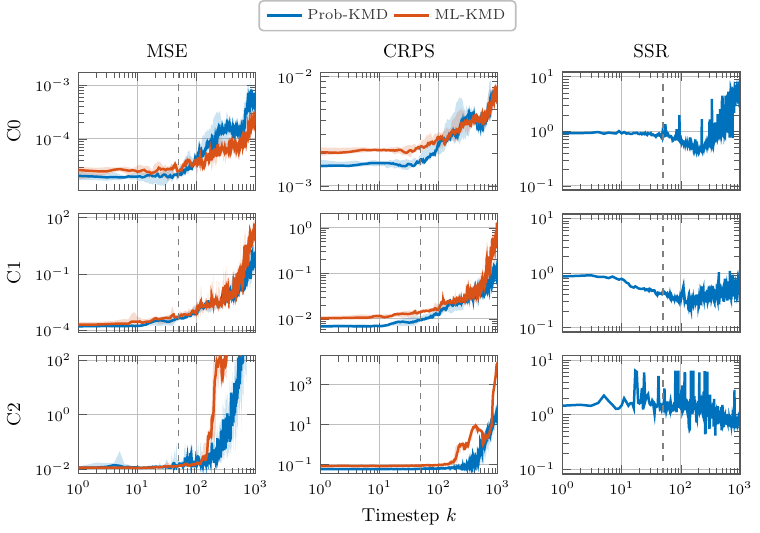}
    \caption{Median and quartile \bk{MSE} and CRPS, and \bk{SSR, against step $k$} for the cavity flow experiments C0 \bk{(no noise), C1 ($\sigma=10^{-2}$)}\bk{ and C2 ($\sigma=10^{-1}$)}\bk{. Dashed line: end of training horizon.}}
    \label{fig:cavity.loss}
\end{figure}
\bk{For experiment C0, Prob-KMD has higher MSE than ML-KMD on rollouts longer than the training interval.}
\bk{In C1 and C2, Prob-KMD has lower MSE than the ML-KMD approximation.}
\bk{However, at large timestep $k$ there is an MSE blow-up in both approximations, particularly in experiment C2 due to eigenvalues with moduli greater than one.}
\bk{Prob-KMD has lower CRPS than ML-KMD for experiments C1 and C2, and at most steps in experiment C0.}

\bk{For C0, the spread--skill ratio is close to 1 for rollouts up to approximately $k=50$.}
\bk{Beyond this, it becomes under-dispersed, before becoming over-dispersed.}
\bk{Under-dispersion starts earlier in experiment C1, even before the length of the training rollout.}
\bk{In experiment C2, the SSR is less smooth, and the model is mildly over-dispersed up to $k\approx500$, with some large spikes (up to $6.4$), and under-dispersed thereafter.}
\bk{Given the smoothness of the CRPS, the spikes are possibly due to occasional samples producing long tails in the state forecast.}

\subsection{\bk{Synthesis}}
\bk{The results illustrate the necessity of probabilistic forecasting.}
\bk{Notably, in the Duffing oscillator D3 experiment, for a trajectory initialised at the unstable fixed point the forecast predicts a multi-modal distribution reflecting the two stable fixed points --- a deterministic model such as the ML-KMD approximation cannot capture this.}

\bk{Prob-KMD attains lower CRPS than the deterministic ML-KMD in all experiments over the training horizon, and generally better over the full predictive horizon except in a few experiments such as V2 whose learnt eigenvalues are erroneously dissipative.}
\bk{The forecasts capture the state data noise but are typically under-dispersed when trained on noisy datasets.}
\bk{The most skilful models may demonstrate over-dispersion, possibly as a consequence of residual injected noise.}

\bk{Numerical results on a non-linear ODE problem with known analytic Koopman eigenfunctions demonstrate that the method is able to identify known eigenfunctions and eigenvalues of the Koopman operator, and that this allows long temporal rollout with accurate forecasts.}
\bk{The two fluid dynamics examples also demonstrate the extension of the method to high-dimensional state spaces (i.e., discretisations of functions).}
\bk{Similarly, the interpretable nature of the model allows verification that the Prob-KMD model can capture the appropriate linear behaviours, e.g., the vortex shedding frequency and the cavity flow lattice.} 

\bk{A comparison of Prob-KMD to a deep ensemble demonstrates a generally comparable performance over the training horizon, but importantly at a significantly reduced computational cost (in this case a factor of the deep ensemble size $\nEnsemble=5$).}
\bk{In experiment A2 with high data noise, the deep ensemble fails to correctly identify part of the Koopman spectrum.}
\bk{Additionally, in experiment D3, DE-KMD does not issue a forecast resolving the distinct basins of attraction for an initial condition at the unstable fixed point of the Duffing oscillator.}
\bk{Deep ensemble approximations do generally give better calibration on long temporal rollouts in a number of experiments.}

\section{Conclusion}
\label{sec:conclusion}
In this work we have outlined the need for calibrated probabilistic forecasting and a method for integrating this into a learnt Koopman mode decomposition (KMD) structure.
A probabilistic observation--linear dynamics--non-linear reconstruction model is proposed, and optimised via the minimisation of the Continuous Ranked Probability Score (CRPS).
Score-based training targets an approximation that samples a distribution that reflects the `true belief' of the forecast.
The method is mostly architecture agnostic: the observation and reconstruction networks in the experiments are chosen to be either multi-layer perceptrons (MLPs) or convolutional neural networks (CNNs) with noise injected as an input.

\bk{Prob-KMD is effective when the linearised approximation of the dynamics is a useful inductive bias, for example, dynamics with decaying and (quasi-)periodic behaviour as exhibited in the examples selected here.}
\bk{The trained Prob-KMD models can identify known Koopman eigenfunctions and spectral properties.}
\bk{Prob-KMD improves on deterministic modelling by offering probabilistic forecasts with lower CRPS, and over the training horizon Prob-KMD gives comparable results to a deep ensemble at reduced computational cost.}
\bk{The Prob-KMD architecture is also able to capture multi-modal distributions, whilst keeping interpretable structured dynamics.}

\bk{Fitting linear dynamics can cause problems on long rollouts for non-dissipative systems as learnt eigenvalues are not constrained to magnitude one.}
\bk{For all models at extremely long rollouts, even small errors in the eigenvalue estimation become amplified.}
Additionally, CRPS only considers marginal distributions for each \bk{state component} and \bk{does not consider} joint structure.
This is partially mitigated by the shared low-dimensional noise input.
\bk{The calibration results are mixed, with Prob-KMD often under-dispersed.}
\bk{Further} evaluation of ensemble forecasting, e.g., rank histograms, is not \bk{included and no direct comparison is made to other similar probabilistic models, e.g.,~\cite{Otto2024,Pan2020} and should be considered in the future.}
{Model training is also expensive, \bk{with each model taking of the order of hours to train}. %
\bk{Training our models is clearly} more \bk{expensive} than the \bk{traditional} numerical DMD and Extended DMD methods, but a cheaper alternative to training a deep ensemble.}

Further work should \bk{investigate} alternative scoring rules to ensure \bk{covariance} is captured, \bk{and the effects of constraining the dynamics on long rollouts}.
In the future, the extension to more complex architectures (e.g., diffusion models, transformers) and more complex test problems is of interest.

\section*{Code and data availability}
The code used to train the models and postprocess the results is provided at \url{github.com/benmkent/prob-kmd}.
The code used to create the datasets for the analytic ODE and Duffing oscillator is provided in the repository.
The vortex shedding dataset is associated with the code repositories for~\cite{Colbrook2024,Colbrook2023res}.
The cavity flow dataset is associated with the code repositories for~\cite{Arbabi2017}. 

\appendix
\section{Model Architecture}\label{sec:architecture}
The model architectures are outlined in \Cref{sec:appendix.mlp,sec:appendix.cnn} with the parameterisation given in \Cref{tab:results.architectures} for each experiment.
\begin{table}
\footnotesize
\centering
\caption{Model architectures for all experiments. ODE experiments (A, D) use MLP architecture and PDE experiments (V, C) use CNN architecture.}
\label{tab:results.architectures}
\setlength{\tabcolsep}{4pt}
\begin{tabular}{llcccc}
\toprule
 & Parameter & A & D & V & C \\
\midrule
\multirow{5}{*}{\rotatebox[origin=c]{90}{Encoder}}
 & Input dim        & $2$              & $2$              & $25\!\times\!100$ & $65\!\times\!65$ \\
 & Output dim       & $4$              & $16$             & $2$              & $16$             \\
 & Noise dim \bk{(Prob-KMD only)}        & $2$              & $4$              & $2$               & $4$              \\
 & MLP hidden dim or CNN channels & $64$             & $64$             & $[16, 32]$        & $[16, 32, 64]$   \\
 & Hidden layers    & $4$              & $5$              & $2$               & $3$               \\
\midrule
\multirow{6}{*}{\rotatebox[origin=c]{90}{Dynamics}}
 & Input dim        & $4$              & $16$             & $2$              & $16$             \\
 & Output dim       & $4$              & $16$             & $2$              & $16$             \\
 & Noise dim \bk{(Prob-KMD only)}       & $2$              & $4$              & $2$               & $4$              \\
 & Hidden dim       & $32$             & $16$             & $6$               & $6$              \\
 & Hidden layers    & $2$              & $3$              & $2$               & $2$              \\
 & Time scaling $\Delta t$       & $0.02$           & $0.02$           & $0.02$            & \bk{$0.02$}            \\
\midrule
\multirow{5}{*}{\rotatebox[origin=c]{90}{Decoder}}
 & Input dim        & $4$              & $16$             & $2$              & $16$             \\
 & Output dim       & $2$              & $2$              & $25\!\times\!100$ & $65\!\times\!65$ \\
 & Noise dim \bk{(Prob-KMD only)}        & $2$              & $4$              & $2$               & $4$              \\
 & MLP hidden dim or CNN channels & $64$             & $64$               & $[128,64,32,16]$  & $[128,64,32,16]$ \\
 & Hidden layers    & $4$              & $5$              & $4$               & $4$               \\
\midrule
& Number of parameters & $28\,682$ & $41\,762$ & $137\,017$ & $248\,833$\\
& Batch size $\nBatch$ \bk{(*0, D3)} & $1000$ & $250$ & $12$ & $12$\\
& Batch size $\nBatch$ \bk{(*1, *2)} & $250$ & $250$ & $10$ & $10$\\
\bottomrule
\end{tabular}
\end{table}

\subsection{Multi-layer Perceptron (MLP)}\label{sec:appendix.mlp}

The observation architecture is a feedforward multi-layer perceptron.
The architecture is defined by
\begin{itemize}
    \item input dimension $\nX$,
    \item output dimension $\nZ$,
    \item noise dimension $\nNoise_\Dict$,
    \item hidden dimension $\nHidden$,
    \item and number of hidden layers $\nLayers$.
\end{itemize}
For each layer $\ell\in\{1,2,\dots,\nLayers-1\}$, the weights $W_\ell \in \R^{\nHidden \times \nHidden}$ and biases $b_\ell \in \R^{\nHidden}$ are trainable parameters.
The input layer is of shape $W_0 \in \R^{\nHidden \times (\nX + \nNoise_\Dict)}$ and $b_0 \in \R^{\nHidden}$, and the output layer is of shape $W_\nLayers \in \R^{\nZ \times \nHidden}$ and $b_\nLayers \in \R^{\nZ}$.
The activation function $\sigma:\R^\nHidden\to\R^\nHidden$ is chosen to be the $\SiLU$ function applied component-wise, where $\SiLU(x) := x(1 + \exp(-x))^{-1}$.
The layers are combined to form the observable function $\Dict_{NN}:\R^{\nX} \times \R^{\nNoise_\Dict} \to \R^\nZ$,
\begin{equation}\label{eq:architecture.mlp}
    \Dict_{NN}(\x, \noise_\Dict) = W_\nLayers \sigma(W_{\nLayers-1} \sigma( \cdots \sigma(W_0 [\x;\, \noise_\Dict] + b_0) \cdots) + b_{\nLayers-1}) + b_\nLayers.
\end{equation}

The same architecture as \eqref{eq:architecture.mlp} is used for the reconstruction function $g_{NN}:\R^{\nZ} \times \R^{\nNoise_g} \to \R^\nX$ with input dimension $\nZ$ and output dimension $\nX$.

\subsection{Convolutional Neural Network (CNN)}\label{sec:appendix.cnn}

An alternative architecture for spatially structured data is a convolutional neural network \bk{(CNN), see e.g.,~\cite{Bishop2023}.}
We consider state spaces $\X=\R^{\nChannels \times \nX_1 \times \nX_2}$ where $\nChannels=1$ is the number of channels for all our test problems, and $\nX_1$ and $\nX_2$ are the spatial dimensions.
The observation architecture is defined by
\begin{itemize}
    \item input dimensions $\nChannels \times \nX_1 \times \nX_2$,
    \item output dimension $\nZ$,
    \item noise dimension $\nNoise_\Dict$,
    \item number of convolutional blocks $\nLayers$,
    \item and number of channels for each block $\nChannels_{i}$ for $i=1,2,\dots,\nLayers$.
\end{itemize}

For an observation, a noise sample is mapped to a $16$-dimensional vector via a trainable affine transformation and the $\SiLU$ activation function, $e := \SiLU(W_e \noise_\Dict + b_e)$ where $W_e \in \R^{\nNoise \times \nNoise_\Dict}, b_e \in \R^{\nNoise}$.
Each convolutional layer $i=1,2,\dots,\nLayers$ applies $\nChannels_i$ convolutional \bk{filters,} \bk{$W_{i,j} \in {\R^{\nChannels_{i-1}\times3\times 3}}$} for $j=1,2,\dots,\nChannels_i$ with a stride length $2$ reducing the spatial extent by a factor of $2$ in each dimension.
\bk{Group normalisation is applied over groups of $8$ channels.}
\bk{Then} a trainable affine transformation parameterised by $\gamma_i:\R^{\nNoise}\to\R^{\nChannels_i},\beta_i:\R^{\nNoise} \to \R^{\nChannels_i}$, conditioned on the noise $e$, is applied for each layer and the $\SiLU$ activation function is applied to each channel output as feature-wise linear modulation (FiLM)~\cite{Perez2018,Alet2025},
\begin{equation}\label{eq:architecture.film}
    y \mapsto \SiLU\left( y \odot (1 + \gamma_i(e)) + \beta_i(e) \right)
\end{equation}
where $\odot$ is component-wise application.
The affine transformation is initialised with $\gamma_i \equiv 0,\beta_i \equiv 0$, i.e., with no perturbations.
Adaptive average pooling reduces the spatial dimension in the final layer to $\nX_{1,\nLayers}:=6$ by \bk{$\nX_{2,\nLayers}:=25$ for the vortex shedding experiments in \Cref{sec:results.vortex} and $\nX_{1,\nLayers}=\nX_{2,\nLayers}=8$ for the cavity flow experiments in \Cref{sec:results.cavity}.} 
A fully connected layer mapping $\R^{\nChannels_{\nLayers}\times \nX_{1,\nLayers} \times \nX_{2,\nLayers}}\to\R^{\nZ}$ returns observable dictionary $\Dict_{\param_{\Dict}}$.

Reconstruction uses the reverse structure. 
First, a fully connected layer maps the observations to a coarse spatial grid of $\nX_{1,0}^g$ times $\nX_{2,0}^{g}$ cells.
\bk{For the vortex experiment this is $4\times13$, for the cavity flow this is $5\times5$.}
Each layer then upsamples by a factor of $2$ with nearest neighbour interpolation.
The same convolutional neural network and FiLM conditioning \eqref{eq:architecture.film} are applied per \bk{layer, but without group normalisation.}
The final \bk{resize} layer maps to the output dimension $\nChannels \times \nX_1 \times \nX_2$, i.e., the state space structure.

\subsection{\bk{Dynamics}}
The linear dynamics use the block diagonal structure \eqref{eq:background.K.rotation}, with learnable scale and rotation parameters $\{(\mu_i, \theta_i)\}_{i=1}^{\nZ/2}$.
To this end, the MLP architecture \eqref{eq:architecture.mlp} with no state input, only noise, is used with two outputs representing a pair $(\mu_i,\theta_i)$.
Independent MLPs are trained for each pair for $i=1,\dots,\nZ/2$.

\subsection{\bk{Deep Ensemble KMD}}
\bk{The DE-KMD approximation is based upon the standard architecture described in~\cite{Lakshminarayanan2017}.}
\bk{In our context, for each ensemble member this learns deterministic MLP observation functions and deterministic dynamics eigenvalues.}
\bk{Probabilistic forecasting comes from the reconstruction which is an MLP that outputs a mean and standard deviation for each state component.}
\bk{A sample from the deep ensemble is formed by pooling samples from the Gaussian reconstruction head of each ensemble member.}
\bk{For the comparison, the same MLP structure and sizing is used for DE-KMD as the ML-KMD model but with double the number of outputs to produce a mean and standard deviation per state component.}

\section{Training}
The models are trained on \bk{a shared} NVIDIA RTX 6000 Ada Generation GPU.
The training uses the AdamW optimiser in \torch{} with the hyperparameters $\beta_1=0.9$, $\beta_2=0.99$, weight decay $10^{-4}$, stability constant $\text{eps}=10^{-7}$.
The batch sizes are given in \Cref{tab:results.architectures}.

For each stage of the curriculum, as listed in \Cref{tab:results.curricula}, the learning rate follows a cosine annealing schedule, with the initial learning rate $\lr$ and the final learning rate $\lr_{\min}$.
An exponential moving average (EMA) of the model is used with factor $0.999$.
The final approximation returned is the EMA model with the best validation loss in the final curriculum stage.
{
\begin{table}
\footnotesize
\centering
\caption{Curriculum training stages for the experiments in \Cref{sec:results}.
$\lr$ shows $[\lr,\lr_{\min}]$ for use of a cosine annealing scheduler within a stage.
For all stages \bk{regularisation weight  $\lossweight_{\regularisation}=10^{-6}$} for A, \bk{D and $\lossweight_{\regularisation}=10^{-3}$} for V, C.
}
\label{tab:results.curricula}
\begin{tabular}{l l r r r r r}
\toprule
& Stage & $\tau$ & Steps & $\lr$
      & $\lossweight_{\Z}$ & $\nSamples$ \\
\midrule
\multirow{5}{*}{A0--A2}
 & Autoencode & $0$   & $10^{4}$          & $[10^{-2},10^{-3}]$ & $0$        & $2$ \\
 & Short      & $5$   & $10^{4}$          & $10^{-3}$           & $10^{-2}$  & $2$ \\
 & Medium     & $25$  & $10^{4}$          & $10^{-3}$           & $10^{-1}$  & $2$ \\
 & Long       & $50$  & $5\cdot10^{4}$    & $[10^{-3},10^{-4}]$ & $10^{-1}$  & $2$ \\
 & Fine-tune  & $50$  & $2\cdot10^{5}$    & $[10^{-3},10^{-4}]$ & $1$        & $5$ \\
\midrule
\multirow{5}{*}{D0--D2}
 & Autoencode & $0$   & $10^{4}$          & $[10^{-2},10^{-3}]$ & $0$        & $2$ \\
 & Short      & $5$   & $10^{4}$          & $10^{-3}$           & $10^{-2}$  & $2$ \\
 & Medium     & $25$  & $10^{4}$          & $10^{-3}$           & $10^{-1}$  & $2$ \\
 & Long       & $50$  & $5\cdot10^{4}$    & $[10^{-3},10^{-4}]$ & $10^{-1}$  & $2$ \\
 & Fine-tune  & $50$  & $2.5\cdot10^{5}$  & $[10^{-3},10^{-4}]$ & $10$       & $5$ \\
\midrule
\multirow{6}{*}{D3}
 & Autoencode & $0$   & $10^{4}$          & $[10^{-2},10^{-3}]$ & $0$        & $2$ \\
 & Short      & $5$   & $10^{4}$          & $10^{-3}$           & $10^{-2}$  & $2$ \\
 & Medium     & $25$  & $10^{4}$          & $10^{-3}$           & $10^{-1}$  & $2$ \\
 & Long       & $50$  & $5\cdot10^{4}$    & $[10^{-3},10^{-4}]$ & $10^{-1}$  & $2$ \\
 & Extra-long & $300$ & $5\cdot10^{4}$    & $[10^{-3},10^{-4}]$ & $1$        & $2$ \\
 & Fine-tune  & $800$ & $10^{5}$          & $[10^{-3},10^{-4}]$ & $10$       & $5$ \\
\midrule
\multirow{5}{*}{V0--V2}
 & Autoencode & $0$   & $2\cdot10^{4}$    & $[10^{-2},10^{-3}]$ & $0$        & $5$ \\
 & Short      & $5$   & $4\cdot10^{4}$    & $10^{-3}$           & $10^{-1}$  & $5$ \\
 & Medium     & $25$  & $2\cdot10^{4}$    & $10^{-3}$           & $10^{-1}$  & $5$ \\
 & Long       & $50$  & $2\cdot10^{4}$    & $[10^{-3},10^{-4}]$ & $10^{-1}$  & $5$ \\
 & Fine-tune  & $50$  & $2\cdot10^{4}$    & $[10^{-3},10^{-6}]$ & $10$       & $5$ \\
\midrule
\multirow{5}{*}{C0--C2}
 & Autoencode & $0$   & $2\cdot10^{4}$    & $[10^{-2},10^{-3}]$ & $0$        & $5$ \\
 & Short      & $5$   & $6\cdot10^{4}$    & $10^{-3}$           & $10^{-2}$  & $5$ \\
 & Medium     & $25$  & $4\cdot10^{4}$    & $10^{-3}$           & $10^{-1}$  & $5$ \\
 & Long       & $50$  & $4\cdot10^{4}$    & $[10^{-3},10^{-4}]$ & $10^{-1}$  & $5$ \\
 & Fine-tune  & $50$  & $1.8\cdot10^{5}$  & $[10^{-3},10^{-6}]$ & $1$        & $5$ \\
\bottomrule
\end{tabular}
\end{table}
}

\section*{Acknowledgements}
\bmk{The authors used Anthropic's Claude tools for proofreading, paper consistency checks and manuscript review.
Additionally, the tools were used to refine the accompanying code (\url{github.com/benmkent/prob-kmd}) and improve typesetting of plots and diagrams.
The authors take full responsibility for the contents of the manuscript, the code and the numerical results.
}

\bibliographystyle{siam/siamplain}
\bibliography{references}

\begin{thebibliography}{10}

\bibitem{Alet2025}
{\sc F.~Alet, I.~Price, A.~El-Kadi, D.~Masters, S.~Markou, T.~R. Andersson,
  J.~Stott, R.~Lam, M.~Willson, A.~Sanchez-Gonzalez, and P.~Battaglia}, {\em
  Skillful joint probabilistic weather forecasting from marginals}, 2025,
  \url{https://arxiv.org/abs/2506.10772}.

\bibitem{Alkin2024}
{\sc B.~Alkin, A.~F\"{u}rst, S.~Schmid, L.~Gruber, M.~Holzleitner, and
  J.~Brandstetter}, {\em Universal physics transformers: A framework for
  efficiently scaling neural operators}, in Advances in Neural Information
  Processing Systems, A.~Globerson, L.~Mackey, D.~Belgrave, A.~Fan, U.~Paquet,
  J.~Tomczak, and C.~Zhang, eds., vol.~37, Curran Associates, Inc., 2024,
  pp.~25152--25194, \url{https://doi.org/10.52202/079017-0793}.

\bibitem{Allen2025}
{\sc A.~Allen, S.~Markou, W.~Tebbutt, J.~Requeima, W.~P. Bruinsma, T.~R.
  Andersson, M.~Herzog, N.~D. Lane, M.~Chantry, J.~S. Hosking, and R.~E.
  Turner}, {\em End-to-end data-driven weather prediction}, Nature, 641 (2025),
  pp.~1172--1179, \url{https://doi.org/10.1038/s41586-025-08897-0}.

\bibitem{Andreuzzi2023}
{\sc F.~Andreuzzi, N.~Demo, and G.~Rozza}, {\em A dynamic mode decomposition
  extension for the forecasting of parametric dynamical systems}, SIAM Journal
  on Applied Dynamical Systems, 22 (2023), pp.~2432--2458,
  \url{https://doi.org/10.1137/22m1481658}.

\bibitem{Arbabi2017}
{\sc H.~Arbabi and I.~Mezi\'{c}}, {\em Study of dynamics in post-transient
  flows using {K}oopman mode decomposition}, Physical Review Fluids, 2 (2017),
  \url{https://doi.org/10.1103/physrevfluids.2.124402}.

\bibitem{Askham2018}
{\sc T.~Askham and J.~N. Kutz}, {\em Variable projection methods for an
  optimized dynamic mode decomposition}, SIAM Journal on Applied Dynamical
  Systems, 17 (2018), pp.~380--416, \url{https://doi.org/10.1137/m1124176}.

\bibitem{Baddoo2022}
{\sc P.~J. Baddoo, B.~Herrmann, B.~J. McKeon, and S.~L. Brunton}, {\em Kernel
  learning for robust dynamic mode decomposition: linear and nonlinear
  disambiguation optimization}, Proceedings of the Royal Society A:
  Mathematical, Physical and Engineering Sciences, 478 (2022),
  \url{https://doi.org/10.1098/rspa.2021.0830}.

\bibitem{Baddoo2023}
{\sc P.~J. Baddoo, B.~Herrmann, B.~J. McKeon, J.~N. Kutz, and S.~L. Brunton},
  {\em Physics-informed dynamic mode decomposition}, Proceedings of the Royal
  Society A: Mathematical, Physical and Engineering Sciences, 479 (2023),
  \url{https://doi.org/10.1098/rspa.2022.0576}.

\bibitem{Bagheri2013}
{\sc S.~Bagheri}, {\em {K}oopman-mode decomposition of the cylinder wake},
  Journal of Fluid Mechanics, 726 (2013), pp.~596--623,
  \url{https://doi.org/10.1017/jfm.2013.249}.

\bibitem{Berkooz1993}
{\sc G.~Berkooz, P.~Holmes, and J.~L. Lumley}, {\em The proper orthogonal
  decomposition in the analysis of turbulent flows}, Annual Review of Fluid
  Mechanics, 25 (1993), pp.~539--575,
  \url{https://doi.org/10.1146/annurev.fl.25.010193.002543}.

\bibitem{Bishop2023}
{\sc C.~M. Bishop and H.~Bishop}, {\em Convolutional networks}, in Deep
  Learning: Foundations and Concepts, Springer International Publishing, 2024,
  pp.~287--324, \url{https://doi.org/10.1007/978-3-031-45468-4_10}.

\bibitem{Bouchacourt2016}
{\sc D.~Bouchacourt, M.~P. Kumar, and S.~Nowozin}, {\em {DISCO} nets:
  Dissimilarity coefficients networks}, in Advances in Neural Information
  Processing Systems, D.~Lee, M.~Sugiyama, U.~von Luxburg, I.~Guyon, and
  R.~Garnett, eds., vol.~29, Curran Associates, Inc., 2016.

\bibitem{Bramburger2024}
{\sc J.~J. Bramburger}, {\em Data-Driven Methods for Dynamic Systems}, Society
  for Industrial and Applied Mathematics, Jan. 2024,
  \url{https://doi.org/10.1137/1.9781611978162}.

\bibitem{Brunton2022}
{\sc S.~L. Brunton, M.~Budi\v{s}i\'{c}, E.~Kaiser, and J.~N. Kutz}, {\em Modern
  {K}oopman theory for dynamical systems}, SIAM Review, 64 (2022),
  pp.~229--340, \url{https://doi.org/10.1137/21m1401243}.

\bibitem{BruntonKutz2022}
{\sc S.~L. Brunton and J.~N. Kutz}, {\em Data-Driven Science and Engineering:
  Machine Learning, Dynamical Systems, and Control}, Cambridge University
  Press, 2nd~ed., May 2022, \url{https://doi.org/10.1017/9781009089517}.

\bibitem{Brunton2016sindy}
{\sc S.~L. Brunton, J.~L. Proctor, and J.~N. Kutz}, {\em Discovering governing
  equations from data by sparse identification of nonlinear dynamical systems},
  Proceedings of the National Academy of Sciences, 113 (2016), pp.~3932--3937,
  \url{https://doi.org/10.1073/pnas.1517384113}.

\bibitem{Brunton2015}
{\sc S.~L. Brunton, J.~L. Proctor, J.~H. Tu, and J.~N. Kutz}, {\em Compressed
  sensing and dynamic mode decomposition}, Journal of Computational Dynamics, 2
  (2015), pp.~165--191, \url{https://doi.org/10.3934/jcd.2015002}.

\bibitem{Budii2012}
{\sc M.~Budi\v{s}i\'{c}, R.~Mohr, and I.~Mezi\'{c}}, {\em Applied
  {K}oopmanism}, Chaos: An Interdisciplinary Journal of Nonlinear Science, 22
  (2012), \url{https://doi.org/10.1063/1.4772195}.

\bibitem{Chen2018}
{\sc R.~T.~Q. Chen, Y.~Rubanova, J.~Bettencourt, and D.~K. Duvenaud}, {\em
  Neural ordinary differential equations}, in Advances in Neural Information
  Processing Systems, S.~Bengio, H.~Wallach, H.~Larochelle, K.~Grauman,
  N.~Cesa-Bianchi, and R.~Garnett, eds., vol.~31, Curran Associates, Inc.,
  2018.

\bibitem{Colbrook2025}
{\sc M.~Colbrook, Z.~Drma\v{c}, and A.~Horning}, {\em An introductory guide to
  {K}oopman learning}, in Operator Theory, D.~Alpay, F.~Colombo, and
  I.~Sabadini, eds., Springer Basel, Oct. 2025, pp.~1--49,
  \url{https://doi.org/10.1007/978-3-0348-0692-3_126-1}.

\bibitem{Colbrook2024}
{\sc M.~J. Colbrook}, {\em The multiverse of dynamic mode decomposition
  algorithms}, in Numerical Analysis Meets Machine Learning, S.~Mishra and
  A.~Townsend, eds., vol.~25 of Handbook of Numerical Analysis, Elsevier, 2024,
  pp.~127--230, \url{https://doi.org/10.1016/bs.hna.2024.05.004}.

\bibitem{Colbrook2023res}
{\sc M.~J. Colbrook, L.~J. Ayton, and M.~Sz\H{o}ke}, {\em Residual dynamic mode
  decomposition: robust and verified {K}oopmanism}, Journal of Fluid Mechanics,
  955 (2023), \url{https://doi.org/10.1017/jfm.2022.1052}.

\bibitem{Colbrook2025rigged}
{\sc M.~J. Colbrook, C.~Drysdale, and A.~Horning}, {\em Rigged dynamic mode
  decomposition: Data-driven generalized eigenfunction decompositions for
  {K}oopman operators}, SIAM Journal on Applied Dynamical Systems, 24 (2025),
  pp.~1150--1190, \url{https://doi.org/10.1137/24m1662370}.

\bibitem{colbrook2024limitspowerskoopmanlearning}
{\sc M.~J. Colbrook, I.~Mezi\'{c}, and A.~Stepanenko}, {\em Limits and powers
  of {K}oopman learning}, 2024, \url{https://arxiv.org/abs/2407.06312}.

\bibitem{Colbrook2023}
{\sc M.~J. Colbrook and A.~Townsend}, {\em Rigorous data-driven computation of
  spectral properties of {K}oopman operators for dynamical systems},
  Communications on Pure and Applied Mathematics, 77 (2024), pp.~221--283,
  \url{https://doi.org/10.1002/cpa.22125}.

\bibitem{Demo2018}
{\sc N.~Demo, M.~Tezzele, and G.~Rozza}, {\em {PyDMD}: Python dynamic mode
  decomposition}, Journal of Open Source Software, 3 (2018),
  \url{https://doi.org/10.21105/joss.00530}.

\bibitem{diaconu2026probabilisticretrofittinglearnedsimulators}
{\sc C.~Diaconu, M.~Cranmer, R.~E. Turner, T.~Marwah, and P.~Mukhopadhyay},
  {\em Probabilistic retrofitting of learned simulators}, in Forty-third
  International Conference on Machine Learning, 2026,
  \url{https://openreview.net/forum?id=nE5g2b8WW5}.

\bibitem{Erichson2019}
{\sc N.~B. Erichson, M.~Muehlebach, and M.~W. Mahoney}, {\em Physics-informed
  autoencoders for {L}yapunov-stable fluid flow prediction}, 2019,
  \url{https://arxiv.org/abs/1905.10866}.

\bibitem{Fortin2014}
{\sc V.~Fortin, M.~Abaza, F.~Anctil, and R.~Turcotte}, {\em Why should ensemble
  spread match the {RMSE} of the ensemble mean?}, Journal of Hydrometeorology,
  15 (2014), pp.~1708--1713, \url{https://doi.org/10.1175/jhm-d-14-0008.1}.

\bibitem{Gal2016}
{\sc Y.~Gal and Z.~Ghahramani}, {\em Dropout as a {B}ayesian approximation:
  representing model uncertainty in deep learning}, in Proceedings of the 33rd
  International Conference on Machine Learning, M.~F. Balcan and K.~Q.
  Weinberger, eds., vol.~48 of Proceedings of Machine Learning Research, PMLR,
  2016, pp.~1050--1059, \url{https://proceedings.mlr.press/v48/gal16.html}.

\bibitem{Garnelo2018}
{\sc M.~Garnelo, D.~Rosenbaum, C.~Maddison, T.~Ramalho, D.~Saxton, M.~Shanahan,
  Y.~W. Teh, D.~Rezende, and S.~M.~A. Eslami}, {\em Conditional neural
  processes}, in Proceedings of the 35th International Conference on Machine
  Learning, J.~Dy and A.~Krause, eds., vol.~80 of Proceedings of Machine
  Learning Research, PMLR, 10--15 Jul 2018, pp.~1704--1713,
  \url{https://proceedings.mlr.press/v80/garnelo18a.html}.

\bibitem{Gneiting2007}
{\sc T.~Gneiting and A.~E. Raftery}, {\em Strictly proper scoring rules,
  prediction, and estimation}, Journal of the American Statistical Association,
  102 (2007), pp.~359--378, \url{https://doi.org/10.1198/016214506000001437}.

\bibitem{Guo2025}
{\sc Y.~Guo, M.~Korda, I.~G. Kevrekidis, and Q.~Li}, {\em Learning parametric
  {K}oopman decompositions for prediction and control}, SIAM Journal on Applied
  Dynamical Systems, 24 (2025), pp.~744--781,
  \url{https://doi.org/10.1137/23m1604576}.

\bibitem{Hersbach2000}
{\sc H.~Hersbach}, {\em Decomposition of the continuous ranked probability
  score for ensemble prediction systems}, Weather and Forecasting, 15 (2000),
  pp.~559--570,
  \url{https://doi.org/10.1175/1520-0434(2000)015<0559:dotcrp>2.0.co;2}.

\bibitem{Ho2020}
{\sc J.~Ho, A.~Jain, and P.~Abbeel}, {\em Denoising diffusion probabilistic
  models}, in Advances in Neural Information Processing Systems, H.~Larochelle,
  M.~Ranzato, R.~Hadsell, M.~F. Balcan, and H.~Lin, eds., vol.~33, Curran
  Associates, Inc., 2020, pp.~6840--6851.

\bibitem{Huhn2023}
{\sc Q.~A. Huhn, M.~E. Tano, J.~C. Ragusa, and Y.~Choi}, {\em Parametric
  dynamic mode decomposition for reduced order modeling}, Journal of
  Computational Physics, 475 (2023),
  \url{https://doi.org/10.1016/j.jcp.2022.111852}.

\bibitem{Koopman1931}
{\sc B.~O. Koopman}, {\em {H}amiltonian systems and transformation in {H}ilbert
  space}, Proceedings of the National Academy of Sciences, 17 (1931),
  pp.~315--318, \url{https://doi.org/10.1073/pnas.17.5.315}.

\bibitem{Koopman1932}
{\sc B.~O. Koopman and J.~von Neumann}, {\em Dynamical systems of continuous
  spectra}, Proceedings of the National Academy of Sciences, 18 (1932),
  pp.~255--263, \url{https://doi.org/10.1073/pnas.18.3.255}.

\bibitem{Korda2018}
{\sc M.~Korda and I.~Mezi\'{c}}, {\em Linear predictors for nonlinear dynamical
  systems: {K}oopman operator meets model predictive control}, Automatica, 93
  (2018), pp.~149--160, \url{https://doi.org/10.1016/j.automatica.2018.03.046}.

\bibitem{Korda2017}
{\sc M.~Korda and I.~Mezi\'{c}}, {\em On convergence of extended dynamic mode
  decomposition to the {K}oopman operator}, Journal of Nonlinear Science, 28
  (2018), pp.~687--710, \url{https://doi.org/10.1007/s00332-017-9423-0}.

\bibitem{Kovachki2023}
{\sc N.~Kovachki, Z.~Li, B.~Liu, K.~Azizzadenesheli, K.~Bhattacharya,
  A.~Stuart, and A.~Anandkumar}, {\em Neural operator: learning maps between
  function spaces with applications to {PDE}s}, Journal of Machine Learning
  Research, 24 (2023), pp.~1--97,
  \url{http://jmlr.org/papers/v24/21-1524.html}.

\bibitem{Lakshminarayanan2017}
{\sc B.~Lakshminarayanan, A.~Pritzel, and C.~Blundell}, {\em Simple and
  scalable predictive uncertainty estimation using deep ensembles}, in Advances
  in Neural Information Processing Systems, I.~Guyon, U.~von Luxburg,
  S.~Bengio, H.~Wallach, R.~Fergus, S.~Vishwanathan, and R.~Garnett, eds.,
  vol.~30, Curran Associates, Inc., 2017.

\bibitem{Lang2026}
{\sc S.~Lang, M.~Alexe, M.~C.~A. Clare, C.~Roberts, R.~Adewoyin,
  Z.~Ben~Bouall\`{e}gue, M.~Chantry, J.~Dramsch, P.~D. Dueben, S.~Hahner,
  P.~Maciel, A.~Prieto-Nemesio, C.~O'Brien, F.~Pinault, J.~Polster, B.~Raoult,
  S.~Tietsche, and M.~Leutbecher}, {\em {AIFS-CRPS}: ensemble forecasting using
  a model trained with a loss function based on the continuous ranked
  probability score}, npj Artificial Intelligence, 2 (2026),
  \url{https://doi.org/10.1038/s44387-026-00073-7}.

\bibitem{Li2017}
{\sc Q.~Li, F.~Dietrich, E.~M. Bollt, and I.~G. Kevrekidis}, {\em Extended
  dynamic mode decomposition with dictionary learning: A data-driven adaptive
  spectral decomposition of the {K}oopman operator}, Chaos: An
  Interdisciplinary Journal of Nonlinear Science, 27 (2017),
  \url{https://doi.org/10.1063/1.4993854}.

\bibitem{Li2020fno}
{\sc Z.~Li, N.~B. Kovachki, K.~Azizzadenesheli, B.~Liu, K.~Bhattacharya,
  A.~Stuart, and A.~Anandkumar}, {\em {F}ourier neural operator for parametric
  partial differential equations}, in International Conference on Learning
  Representations, 2021, \url{https://openreview.net/forum?id=c8P9NQVtmnO}.

\bibitem{liu2024}
{\sc T.~J. Liu, N.~Boull\'{e}, R.~Sarfati, and C.~Earls}, {\em {LLM}s learn
  governing principles of dynamical systems, revealing an in-context neural
  scaling law}, in Proceedings of the 2024 Conference on Empirical Methods in
  Natural Language Processing, Y.~Al-Onaizan, M.~Bansal, and Y.-N. Chen, eds.,
  Miami, Florida, USA, Nov. 2024, Association for Computational Linguistics,
  pp.~15097--15117, \url{https://doi.org/10.18653/v1/2024.emnlp-main.842}.

\bibitem{Lu2021}
{\sc L.~Lu, P.~Jin, G.~Pang, Z.~Zhang, and G.~E. Karniadakis}, {\em Learning
  nonlinear operators via {DeepONet} based on the universal approximation
  theorem of operators}, Nature Machine Intelligence, 3 (2021), pp.~218--229,
  \url{https://doi.org/10.1038/s42256-021-00302-5}.

\bibitem{Lusch2018}
{\sc B.~Lusch, J.~N. Kutz, and S.~L. Brunton}, {\em Deep learning for universal
  linear embeddings of nonlinear dynamics}, Nature Communications, 9 (2018),
  \url{https://doi.org/10.1038/s41467-018-07210-0}.

\bibitem{Mezic2005}
{\sc I.~Mezi\'{c}}, {\em Spectral properties of dynamical systems, model
  reduction and decompositions}, Nonlinear Dynamics, 41 (2005), pp.~309--325,
  \url{https://doi.org/10.1007/s11071-005-2824-x}.

\bibitem{Mishra2024}
{\sc S.~Mishra and A.~Townsend}, eds., {\em Numerical Analysis Meets Machine
  Learning}, vol.~25 of Handbook of Numerical Analysis, Elsevier, 2024.

\bibitem{Morton2019}
{\sc J.~Morton, F.~D. Witherden, and M.~J. Kochenderfer}, {\em Deep variational
  {K}oopman models: inferring {K}oopman observations for uncertainty-aware
  dynamics modeling and control}, in Proceedings of the Twenty-Eighth
  International Joint Conference on Artificial Intelligence, IJCAI-19,
  International Joint Conferences on Artificial Intelligence Organization,
  2019, pp.~3173--3179, \url{https://doi.org/10.24963/ijcai.2019/440}.

\bibitem{Otto2024}
{\sc S.~Otto, S.~Peitz, and C.~Rowley}, {\em Learning bilinear models of
  actuated {K}oopman generators from partially observed trajectories}, SIAM
  Journal on Applied Dynamical Systems, 23 (2024), pp.~885--923,
  \url{https://doi.org/10.1137/22m1523601}.

\bibitem{Otto2019}
{\sc S.~E. Otto and C.~W. Rowley}, {\em Linearly recurrent autoencoder networks
  for learning dynamics}, SIAM Journal on Applied Dynamical Systems, 18 (2019),
  pp.~558--593, \url{https://doi.org/10.1137/18m1177846}.

\bibitem{Pan2020}
{\sc S.~Pan and K.~Duraisamy}, {\em Physics-informed probabilistic learning of
  linear embeddings of nonlinear dynamics with guaranteed stability}, SIAM
  Journal on Applied Dynamical Systems, 19 (2020), pp.~480--509,
  \url{https://doi.org/10.1137/19m1267246}.

\bibitem{Perez2018}
{\sc E.~Perez, F.~Strub, H.~de~Vries, V.~Dumoulin, and A.~Courville}, {\em
  {FiLM}: Visual reasoning with a general conditioning layer}, Proceedings of
  the {AAAI} Conference on Artificial Intelligence, 32 (2018),
  \url{https://doi.org/10.1609/aaai.v32i1.11671}.

\bibitem{Price2025}
{\sc I.~Price, A.~Sanchez-Gonzalez, F.~Alet, T.~R. Andersson, A.~El-Kadi,
  D.~Masters, T.~Ewalds, J.~Stott, S.~Mohamed, P.~Battaglia, R.~Lam, and
  M.~Willson}, {\em Probabilistic weather forecasting with machine learning},
  Nature, 637 (2025), pp.~84--90,
  \url{https://doi.org/10.1038/s41586-024-08252-9}.

\bibitem{Proctor2016}
{\sc J.~L. Proctor, S.~L. Brunton, and J.~N. Kutz}, {\em Dynamic mode
  decomposition with control}, SIAM Journal on Applied Dynamical Systems, 15
  (2016), pp.~142--161, \url{https://doi.org/10.1137/15m1013857}.

\bibitem{Raissi2018}
{\sc M.~Raissi, Z.~Wang, M.~S. Triantafyllou, and G.~E. Karniadakis}, {\em Deep
  learning of vortex-induced vibrations}, Journal of Fluid Mechanics, 861
  (2019), pp.~119--137, \url{https://doi.org/10.1017/jfm.2018.872}.

\bibitem{Rowley2009}
{\sc C.~W. Rowley, I.~Mezi\'{c}, S.~Bagheri, P.~Schlatter, and D.~S.
  Henningson}, {\em Spectral analysis of nonlinear flows}, Journal of Fluid
  Mechanics, 641 (2009), pp.~115--127,
  \url{https://doi.org/10.1017/s0022112009992059}.

\bibitem{Sashidhar2022}
{\sc D.~Sashidhar and J.~N. Kutz}, {\em Bagging, optimized dynamic mode
  decomposition for robust, stable forecasting with spatial and temporal
  uncertainty quantification}, Philosophical Transactions of the Royal Society
  A: Mathematical, Physical and Engineering Sciences, 380 (2022),
  \url{https://doi.org/10.1098/rsta.2021.0199}.

\bibitem{Schmid2010}
{\sc P.~J. Schmid}, {\em Dynamic mode decomposition of numerical and
  experimental data}, Journal of Fluid Mechanics, 656 (2010), pp.~5--28,
  \url{https://doi.org/10.1017/s0022112010001217}.

\bibitem{Scholz2026}
{\sc J.~Scholz, C.~Diaconu, A.~Shysheya, S.~Markou, and R.~E. Turner}, {\em
  {OtterWeather}: Highly skillful medium-range weather forecasting on a single
  {GPU}}, in AI{\&}PDE: ICLR 2026 Workshop on AI and Partial Differential
  Equations, 2026, \url{https://openreview.net/forum?id=B2FtMotBw4}.

\bibitem{Schrab2023}
{\sc A.~Schrab, I.~Kim, M.~Albert, B.~Laurent, B.~Guedj, and A.~Gretton}, {\em
  {MMD} aggregated two-sample test}, Journal of Machine Learning Research, 24
  (2023), pp.~1--81, \url{http://jmlr.org/papers/v24/21-1289.html}.

\bibitem{Shankar2000}
{\sc P.~N. Shankar and M.~D. Deshpande}, {\em Fluid mechanics in the driven
  cavity}, Annual Review of Fluid Mechanics, 32 (2000), pp.~93--136,
  \url{https://doi.org/10.1146/annurev.fluid.32.1.93}.

\bibitem{Shen1991}
{\sc J.~Shen}, {\em Hopf bifurcation of the unsteady regularized driven cavity
  flow}, Journal of Computational Physics, 95 (1991), pp.~228--245,
  \url{https://doi.org/10.1016/0021-9991(91)90261-i}.

\bibitem{Takeishi2017}
{\sc N.~Takeishi, Y.~Kawahara, Y.~Tabei, and T.~Yairi}, {\em {B}ayesian dynamic
  mode decomposition}, in Proceedings of the Twenty-Sixth International Joint
  Conference on Artificial Intelligence, IJCAI-17, International Joint
  Conferences on Artificial Intelligence Organization, Aug. 2017,
  pp.~2814--2821, \url{https://doi.org/10.24963/ijcai.2017/392}.

\bibitem{Takeishi2017b}
{\sc N.~Takeishi, Y.~Kawahara, and T.~Yairi}, {\em Learning {K}oopman invariant
  subspaces for dynamic mode decomposition}, in Advances in Neural Information
  Processing Systems, I.~Guyon, U.~von Luxburg, S.~Bengio, H.~Wallach,
  R.~Fergus, S.~Vishwanathan, and R.~Garnett, eds., vol.~30, Curran Associates,
  Inc., 2017.

\bibitem{Tipping1999}
{\sc M.~E. Tipping and C.~M. Bishop}, {\em Probabilistic principal component
  analysis}, Journal of the Royal Statistical Society Series B: Statistical
  Methodology, 61 (1999), pp.~611--622,
  \url{https://doi.org/10.1111/1467-9868.00196}.

\bibitem{Trefethen2000}
{\sc L.~N. Trefethen}, {\em Spectral Methods in {MATLAB}}, Society for
  Industrial and Applied Mathematics, 2000,
  \url{https://doi.org/10.1137/1.9780898719598}.

\bibitem{Tu2014}
{\sc J.~H. Tu, C.~W. Rowley, D.~M. Luchtenburg, S.~L. Brunton, and J.~N. Kutz},
  {\em On dynamic mode decomposition: Theory and applications}, Journal of
  Computational Dynamics, 1 (2014), pp.~391--421,
  \url{https://doi.org/10.3934/jcd.2014.1.391}.

\bibitem{Turner2011}
{\sc R.~E. Turner and M.~Sahani}, {\em Two problems with variational
  expectation maximisation for time series models}, in Bayesian Time Series
  Models, D.~Barber, A.~T. Cemgil, and S.~Chiappa, eds., Cambridge University
  Press, Aug. 2011, pp.~104--124,
  \url{https://doi.org/10.1017/cbo9780511984679.006}.

\bibitem{Williams2015}
{\sc M.~O. Williams, I.~G. Kevrekidis, and C.~W. Rowley}, {\em A data-driven
  approximation of the {K}oopman operator: Extending dynamic mode
  decomposition}, Journal of Nonlinear Science, 25 (2015), pp.~1307--1346,
  \url{https://doi.org/10.1007/s00332-015-9258-5}.

\bibitem{Williams2015b}
{\sc M.~O. Williams, C.~W. Rowley, and I.~G. Kevrekidis}, {\em A kernel-based
  method for data-driven {K}oopman spectral analysis}, Journal of Computational
  Dynamics, 2 (2015), pp.~247--265, \url{https://doi.org/10.3934/jcd.2015005}.

\bibitem{Yeung2019}
{\sc E.~Yeung, S.~Kundu, and N.~Hodas}, {\em Learning deep neural network
  representations for {K}oopman operators of nonlinear dynamical systems}, in
  2019 American Control Conference (ACC), 2019, pp.~4832--4839,
  \url{https://doi.org/10.23919/ACC.2019.8815339}.

\bibitem{Zamo2018}
{\sc M.~Zamo and P.~Naveau}, {\em Estimation of the continuous ranked
  probability score with limited information and applications to ensemble
  weather forecasts}, Mathematical Geosciences, 50 (2018), pp.~209--234,
  \url{https://doi.org/10.1007/s11004-017-9709-7}.

\end{thebibliography}

\end{document}